\documentclass[a4paper]{amsproc}

\usepackage{amsfonts}
\usepackage{amssymb}
\usepackage{amsmath}

 \usepackage[colorlinks=true]{hyperref}
 \hypersetup{pdftex,
            breaklinks,
            bookmarks,
            colorlinks=true,
            citecolor=black,
            linkcolor=black,
            urlcolor=black,
            pagecolor=black,
            filecolor=black,
            pdftitle={The transfer operator for the Hecke triangle groups},
            pdfauthor={D. Mayer, T. Mhlenbruch, F. Strmberg},
            pdfkeywords={transfer operator, Hecke triangle group, continued fraction},
            pdfstartview=FitH,
            unicode,
            plainpages=false}

\newtheorem{theorem}{Theorem}[section]

\newtheorem{lemma}[theorem]{Lemma}
\newtheorem{proposition}{Proposition}

\theoremstyle{definition}

\newtheorem{remark}{Remark}

\newcommand{\IA}{\ensuremath{\mathbb{A}}}
\newcommand{\IC}{\ensuremath{\mathbb{C}}}
\newcommand{\IP}{\ensuremath{\mathbb{P}}}

\newcommand{\IR}{\ensuremath{\mathbb{R}}}
\newcommand{\IZ}{\ensuremath{\mathbb{Z}}}
\newcommand{\IN}{\ensuremath{\mathbb{N}}}
\newcommand{\IH}{\ensuremath{\mathbb{H}}}

\newcommand{\cA}{\ensuremath{\mathcal{A}}}
\newcommand{\cB}{\ensuremath{\mathcal{B}}}

\newcommand{\cF}{\ensuremath{\mathcal{F}}}

\newcommand{\cM}{\ensuremath{\mathcal{M}}}
\newcommand{\cN}{\ensuremath{\mathcal{N}}}
\newcommand{\cP}{\ensuremath{\mathcal{P}}}

\newcommand{\dd}{\ensuremath{\mathrm{d}}}
\newcommand{\id}{\ensuremath{\mathbf{1}}}

\newcommand{\CF}{{$\lambda_q$-CF}}
\newcommand{\Ar}{\ensuremath{\cA_q^\text{reg}}}
\newcommand{\Adr}{\ensuremath{\cA_q^\text{dreg}}}
\newcommand{\Iq}{\ensuremath{I_q}}
\newcommand{\Ir}{\ensuremath{I_{R_q}}}
\newcommand{\Fix}{\ensuremath{\operatorname{Fix}}}

\newcommand{\lb}{\ensuremath{[\![}} 
\newcommand{\rb}{\ensuremath{]\!]}} 
\newcommand{\lbs}{\ensuremath{[}} 
\newcommand{\rbs}{\ensuremath{]}} 
\newcommand{\lbd}{\ensuremath{[\![}} 
\newcommand{\rbd}{\ensuremath{]\!]^\star}} 

\newcommand{\re}[1]{\ensuremath{{\mathrm{Re}\left(#1\right)}}}

\newcommand{\SL}[1]{\ensuremath{{\mathrm{SL}\!\left(2, #1 \right)}}}

\newcommand{\PSL}[1]{\ensuremath{{\mathrm{PSL}\!\left(2, #1 \right)}}}

\newcommand{\Matrix}[4]{{\begin{pmatrix} #1 & #2 \\ #3 & #4 \end{pmatrix}}}
\newcommand{\Tr}[1]{\ensuremath{{\mathrm{Tr}\left({#1}\right)}}}

\newcommand\txtfrac[2]{{\textstyle \frac{#1}{#2}}}
\newcommand{\abs}[1]{\ensuremath{{\left\lvert#1\right\rvert}}}
\newcommand{\ov}[1]{\ensuremath{\overline{ #1 }}}
\newcommand{\nextinteger}[1]{\ensuremath{\left\lfloor #1\right\rfloor}}
\newcommand{\nli}[1]{\ensuremath{\left\langle #1 \right\rangle_q}}
\newcommand{\nlid}[1]{\ensuremath{\left\langle #1 \right\rangle_q^\star}}

\title[Transfer operator for Hecke triangle groups]
      {The transfer operator for the Hecke triangle groups}

\author[Dieter Mayer, Tobias M\"uhlenbruch and Fredrik Str\"omberg]{}


 \email{dieter.mayer@tu-clausthal.de}
 \email{tobias.muehlenbruch@fernuni-hagen.de}
 \email{stroemberg@mathematik.tu-darmstadt.de}
\thanks{The second author was supported by German Research Foundation (DFG) grant Ma 633/16-1: Transfer operators and non-arithmetic quantum chaos}

\begin{document}

\newlength{\wone}
\settowidth{\wone}{\scriptsize{1}}
\newlength{\wi}
\settowidth{\wi}{\scriptsize{i}}
\newlength{\wonemi}
\setlength{\wonemi}{\wone}
\addtolength{\wonemi}{-\wi}

\maketitle

\centerline{\scshape Dieter Mayer}
\medskip
{\footnotesize
\centerline{Lower Saxony Professorship}
\centerline{Institute for Theoretical Physics, TU Clausthal}
   \centerline{D-38678 Clausthal-Zellerfeld, Germany}
} 

\medskip

\centerline{\scshape Tobias M\"uhlenbruch }
\medskip
{\footnotesize
 \centerline{Department of Mathematics and Computer Science, FernUniversit\"at in Hagen}
   \centerline{D-58084 Hagen, Germany}
} %

\medskip

\centerline{\scshape Fredrik Str\"omberg}
\medskip
{\footnotesize
 \centerline{Department of Mathematics, TU Darmstadt}
   \centerline{D-64289 Darmstadt, Germany}
} %
\bigskip


\begin{abstract}
In this paper we extend the transfer operator approach to Selberg's zeta function for cofinite
Fuchsian groups to the Hecke triangle groups $G_q,\, q=3,4,\ldots$, which are non-arithmetic for
$q\not= 3,4,6$. For this we make use of a Poincar\'e map for the geodesic flow on the corresponding
Hecke surfaces, which has been constructed in \cite{MS08}, and which is closely related to the natural extension of the generating map for the
so-called Hurwitz-Nakada continued fractions.
We also derive functional equations for the eigenfunctions of the transfer operator which for
eigenvalues $\rho =1$ are expected to be closely related to the period functions of Lewis and Zagier
for these Hecke triangle groups.
\end{abstract}
\section{Introduction}
\label{A}
This paper continues the study of the transfer operator for cofinite Fuchsian groups and their Selberg
zeta functions in e.g.~\cite{CM00,CM01}.
By \emph{modular groups} we mean finite index subgroups of the modular group $\PSL{\IZ}$.
For such groups the transfer operator approach to Selberg's zeta function \cite{CM00} has led to interesting new developments in number theory,
 like the theory of period functions for Maa{\ss} wave forms by Lewis and Zagier \cite{LZ01}.
One would like to extend this theory to more general Fuchsian groups, especially the nonarithmetic
ones. One possibility to obtain such a generalization is via a cohomological approach \cite{BLZ09}, which
has recently been considered for the case $G_3 =\PSL{\IZ}$ in \cite{BM09}. We concentrate on
the transfer operator approach to this circle of problems and started to work out this approach in \cite{MM09, MS08} for the Hecke triangle groups which, contrary to modular groups studied up to now, are mostly non-arithmetic.

The transfer operator was introduced by D.~Ruelle \cite{Ru94} to (primarily) investigate analytic properties of dynamical zeta functions.
A typical example of such a function is the Selberg zeta function $Z_S(s)$ for the geodesic flow on a surface
of constant negative curvature, which connects the length spectrum of this flow with spectral properties of the corresponding Laplacian. It is defined by
\begin{equation}
\label{Selberg zeta}
Z_S(s) = \prod_{\gamma} \prod_{k=0}^\infty \left(1-e^{-(s+k)l(\gamma)} \right),
\end{equation}
where the outer product is taken over all prime periodic orbits $\gamma$ of period $l(\gamma)$ of the geodesic flow on the unit tangent bundle of the surface.
The period coincides in this case with the length of the corresponding closed geodesic.
If $\cP:\Sigma \to \Sigma$ is the Poincar\'e map on a section $\Sigma$ of the flow $\Phi_t$ Ruelle showed that $Z_S(s)$ can be rewritten as
\[
Z_S(s) = \prod_{k=0}^\infty \frac{1}{\zeta_R(s+k)},
\]
where $\zeta_R$ denotes the Ruelle zeta function for the Poincar\'e map $\cP$, defined as
\begin{equation*}
\label{Ruelle zeta}
\zeta_R(s) = \exp \left( \sum_{n=1}^\infty \frac{1}{n} Z_n(s) \right)
\end{equation*}
where
\begin{equation*}
\label{Ruelle zeta 2}
Z_n(s) = \sum_{x \in \Fix \cP^n} \exp \left(-s \sum_{k=0}^{n-1} r\left(\cP^k(x)\right) \right), n\ge1,
\end{equation*}
are the so called dynamical partition functions and $r:\Sigma \to \IR_+$ is the recurrence time function with respect to the map $\cP$, defined through
\[
\Phi_{r(x)}(x) \in \Sigma \text{ for } x \in \Sigma
\quad \text{and} \quad
\Phi_t(x) \notin \Sigma \text{ for } 0 < t < r(x).
\]
In the transfer operator approach, the Selberg zeta function is expressed in terms of the Fredholm determinant of an operator $\mathcal{L}_s$ as $Z_S(s)=\det (1-\mathcal{L}_s)$.
From this relation it is clear that the zeros of $Z_S(s)$ are directly related to the values of $s$ for which $\mathcal{L}_s$ has eigenvalue one.
Furthermore, for modular groups, the corresponding eigenfunctions, in a certain Banach space of holomorphic functions,
 can be directly related to certain automorphic functions (cf.~e.g.~\cite{CM01}).
In this paper we work out the details of the transfer operator approach for the Hecke triangle groups, $G_q$,
and the corresponding surfaces $\cM_q$ (these will be defined more precisely in the next section).

A Poincar\'e map and cross-section for the geodesic flow on the Hecke surfaces was constructed in \cite{MS08}.
Additionally, this Poincar\'e map was also shown to be closely related to the natural extension of the generating map $f_q:I_q\to I_q$
for a particular type of continued fraction expansions, denoted by $\lambda_q$-continued fractions or \CF\,'s for short.
In this precise form and generality these were first considered by Nakada \cite{Na95}
but since they are based on the nearest $\lambda_q$ multiple map $f_q$ they also generalize the nearest integer expansions considered by Hurwitz \cite{Hu89}
and possess many of the same characteristics as these. Therefore we sometimes also call them Hurwitz-Nakada continued fractions.
There is also a close relationship to the Rosen $\lambda$-continued fractions \cite{Ro54,RS92,SS95}.
For a precise description of this relationship see \cite[Rem.~15]{MS08}.

For the modular surfaces a Poincar\'e map was constructed through the natural extension of the Gau{\ss} map $T_G:[0,1)\to[0,1),\, T_G(x)=\frac{1}{x}\mod 1,\, x\neq 0,
$ which is related to the so-called simple continued fractions (cf.~e.g.~\cite{Ma90,CM00,CM01}).
In contrast to this case, in the present case (i.e. for Hecke surfaces $\cM_q$), as we will see, there is not a one-to-one
correspondence between the periodic orbits of the map $f_q$ generating the \CF\,'s and the periodic orbits of the geodesic
flow $\Phi_t$. Indeed, for every $G_q$, there exist two periodic points, $r_q$ and $-r_q$, of $f_q$, which correspond to the same periodic orbit
$\mathcal{O}$ for the geodesic flow. For $q=3$, i.e. for the modular group, this fact follows by the results of Hurwitz \cite{Hu89}, who discovered the
existence and properties of these two periodic points.
As a consequence, the Fredholm determinant
of the Ruelle transfer operator, $\mathcal{L}_s$, for the Hurwitz-Nakada map $f_q$ contains the contribution of the closed orbit $\mathcal{O}$ twice.
Therefore it does not, by itself, correctly describe the corresponding Selberg zeta function
(\ref{Selberg zeta}) in same manner as in e.g.~\cite{CM00} for the modular groups.
To correct for this overcounting we introduce another transfer operator, $\mathcal{K}_s$, whose Fredholm
determinant exactly corresponds to the contribution of the orbit $\mathcal{O}$ to $Z_S(s)$.
The form of this operator can be deduced directly from the \CF\, expansion of the point
$r_q$ and furthermore its spectrum can be determined explicitly, leading to
regularly spaced zeros of its Fredholm determinant, $\text{det} (1-\mathcal{K}_s)$, in the complex
$s$-plane. In Section \ref{secK} we will use the operator $\mathcal{K}_s$ to show the following formula for the Selberg zeta function for Hecke triangle groups:
\begin{equation}\label{LoverK}
 Z_S(s) = \frac{\text{det}(1-\mathcal{L}_s)}{\text{det}(1-\mathcal{K}_s)}.
\end{equation}
As in the case of the modular surfaces and the Gau{\ss} map $T_G$, we will see that the holomorphic eigenfunctions
 of the transfer operator $\mathcal{L}_s$ fulfil certain functional equations with a finite number of terms.
In the case $q=3$ it was recently shown \cite{BM09}, that for $0< \re{s}<1, s\not=\frac{1}{2},$
there is a one-to-one correspondence between eigenfunctions of $\mathcal{L}_{s}$ with eigenvalue $1$ (satisfying certain additional conditions) and Maa{\ss} waveforms, i.e. square-integrable
eigenfunctions of the Laplace-Beltrami operator, on the modular surface $\mathcal{M}_{3}$.
This relationship can be interpreted as a correspondence between the classical dynamics, in the guise of the geodesic flow,
and the quantum mechanical dynamics on the surface $\mathcal{M}_3$.
Since the connection between the transfer operator and the geodesic flow is just as strong for any Hecke triangle group, $G_q$, as it is for $G_3$ we expect similar relationships to hold
in general. That is, for any $q \ge 3$ we expect some form of explicit correspondence between holomorphic eigenfunctions of $\mathcal{L}_s$ with eigenvalue $\rho(s) =1$ and automorphic functions of $G_q$.
Observe that almost all Hecke triangle groups, although possessing an infinite number of Maa{\ss} waveforms, due to being cycloidal, are non arithmetic.
Therefore such a correspondence would extend the transfer operator approach to the theory of period functions of Lewis and Zagier
\cite{LZ01} to a whole class of non-arithmetic Fuchsian groups. We hope to come back to this question soon.

An outline of the remaining part of present paper is as follows: In Section~\ref{C} we introduce the Hecke triangle
groups and recall the necessary properties of the $\lambda_q$-continued fractions, including the construction of a Markov partition
for the corresponding generating map $f_q$.
In Section~\ref{D} we recall properties of the geodesic flow on the unit tangent bundle of the Hecke surface $\cM_q$. We also
briefly repeat the results from \cite{MS08} concerning the construction of the corresponding Poincar\'e section $\Sigma$ and the Poincar\'e map $\cP:\Sigma\to \Sigma$.
The definitions of, and relationships between the Ruelle and Selberg zeta functions are also presented.
The transfer operator $\mathcal{L}_s$ for the map $f_q$ is discussed in
Section~\ref{E}. We show that it is a nuclear operator when acting on a certain Banach space $B$ of
vector-valued holomorphic functions, determined by the Markov partitions for $f_q$. We also show that is
has a meromorphic extension to the entire complex $s$-plane.
In Section~\ref{F} we define a symmetry operator, $P:B\to B$, commuting
with the transfer operator. This allows us to restrict the operator $\mathcal{L}_s$ to the two eigenspaces
$B_\epsilon,\,\epsilon =\pm 1 $ of $P$. Using this restriction, $\mathcal{L}_{s,á¸³\epsilon}$, we derive scalar functional equations
fulfilled by the eigenfunctions with eigenvalue $1$.
In Section~\ref{G} we give more details about the Ruelle and Selberg zeta functions.
We also prove the exact relationship between the zeta functions and the corresponding transfer operators $\mathcal{L}_s$ and $\mathcal{K}_s$, arriving finally at \eqref{LoverK} by means of Theorem \ref{main-theorem}.

\begin{flushleft}
\end{flushleft}\section{Background}\label{C}
\subsection{The Hecke triangle groups} \label{C0}
Consider the \emph{projective special linear group}
\[
\PSL{\IR} = \SL{\IR} / \left\{ \pm \id \right\},
\]
where $\SL{\IR}$, the \emph{special linear group}, consists of $2\times 2$ matrices with real entries and determinant $1$ and $\id$ is the $2\times 2$ identity matrix.
For notational convenience, elements of $\PSL{\IR}$ are usually identified with their matrix representatives in $\SL{\IR}$.
If $\IH$ denotes the hyperbolic upper half-plane, that is $\IH=\{z=x+iy \mid y >0\}$ together with the hyperbolic arc length $ds=y^{-1}|dz|$ and area measure $d\mu=y^{-2}dxdy$
then $\PSL{\IR}$ can be identified with the group of orientation preserving isometries of $\IH$. The action of $\PSL{\IR}$ on $\IH$ is given by \emph{M{\"o}bius transformations}
\begin{equation*}
\label{D1.1}
g \, z := \frac{az+b}{cz+d}\quad \text{for} \quad g=\Matrix{a}{b}{c}{d}\in \PSL{\IR}.
\end{equation*}
One can easily verify that this is indeed an orientation preserving (and conformal) action on $\IH$, which additionally preserves the hyperbolic arc length $ds$
as well as the area measure $d\mu$. Furthermore, it extends to an action on the boundary of $\IH$, $\partial \IH = \IP_\IR=\IR\cup\{\infty\}$, the projective real line.

Of particular interest in the theory of automorphic forms are the cofinite Fuchsian groups, that is discrete subgroups $\Gamma$ of $\PSL{\IR}$,
where the corresponding quotient orbifold $\Gamma \backslash \IH$ has finite hyperbolic area.
Although, strictly speaking, these orbifolds are not in general surfaces according to modern terminology (they have marked points),
 we nevertheless view these as Riemann surfaces in the classical sense, i.e. they possess a (not necessarily everywhere smooth) Riemannian structure.

For an integer $q \geq 3$, the \emph{Hecke triangle group}, $G_q$, is the cofinite Fuchsian group generated by
\begin{equation}
\label{C0.3}
S = \Matrix{0}{-1}{1}{0}:z\mapsto -\frac{1}{z}
\quad \text{and} \quad
T:=T_q = \Matrix{1}{\lambda_q}{0}{1}:z\mapsto z+\lambda_q,
\end{equation}
where
\begin{equation*}
\label{C0.1}
\lambda_q = 2 \cos \left( \frac{\pi}{q} \right) \in [1,2).
\end{equation*}
The only relations in the group $G_q$ are given by
\begin{equation*}
\label{C0.2}
S^2 = (ST_q)^q = \id.
\end{equation*}
The orbifold $\cM_q= G_q \backslash \IH$ is usually said to be a \emph{Hecke triangle surface} and
we denote by $\pi:\IH \to \cM_q$ the natural projection map $\pi(z)=G_q z$.
For practical purposes, $\cM_q$ is usually identified with the
standard fundamental domain of $G_q$
\begin{equation*}
\label{D1.3}
\cF_q = \big\{z \in \IH \mid \abs{z} \geq 1, \abs{\re{z}} \leq \txtfrac{\lambda_q}{2} \big\}
\end{equation*}
where the sides are pairwise identified by the generators in (\ref{C0.3}).

We say that two points $x,y\in\IH\cup\IP_\IR$ are \emph{$G_q$-equivalent} if there exists a $g \in G_q$ such that $x = g \, y$.

An element $g=\begin{pmatrix}a&b\cr c&d\cr \end{pmatrix} \in \PSL{\IR}$ is called \emph{elliptic}, \emph{hyperbolic} or \emph{parabolic} depending on whether
$\abs{\Tr{g}} := \abs{a +d} <2$, $>2$ or $= 2$.
The same notation applies for the fixed points of the corresponding M\"obius transformation.
Note that the type of fixed point is preserved under conjugation, $g \mapsto A g A^{-1}$, by $A \in \PSL{\IR}$.
A parabolic fixed point is a degenerate fixed point, belongs to $\IP_\IR$, and is usually called a \emph{cusp}.
Elliptic fixed points appear in pairs, $z$ and $\ov{z}$ with $z \in \IH$.
Hyperbolic fixed points also appear in pairs $x, x^\star \in \IP_\IR$, where $x^\star$ is
said to be the repelling conjugate point of the attractive fixed point $x$.

\subsection{$\lambda_q$-continued fractions}
\label{C1}
Consider finite or infinite sequences $\lbs a_i\rbs_i$ with $a_i\in\mathbb{Z}$ for all $i$.
We denote a periodic subsequence within an infinite sequence by overlining the periodic part and a finitely often repeated pattern is denoted by a power, where the power $^0$ means absence of the pattern, hence
\begin{align*}
\lbs a_1, \ov{a_2, a_3} \rbs
  &= \lbs a_1, \,a_2, a_3,\,a_2, a_3,\,a_2, a_3,\,\ldots \rbs, \\
\lbs a_1, (a_2, a_3)^i,a_4, \ldots \rbs
&= \lbs a_1, \underbrace{a_2, a_3,a_2, a_3, \ldots, a_2,a_3}_{i \, \mathrm{times} \; a_2,a_3}, a_4, \ldots \rbs \quad \text{and} \\
\lbs a_1, (a_2)^0, a_3, \ldots \rbs
  &= \lbs a_1, a_3, \ldots \rbs.\\
\intertext{Furthermore, by the negative of a sequence we mean the following:}
-\lbs a_1,a_2,\ldots \rbs &= \lbs -a_1,-a_2,\ldots \rbs.
\end{align*}

Put
\begin{equation*}
\label{C1.3}
h_q :=
\begin{cases}
\frac{q-2}{2} \quad & \text{for even $q$ and} \\
\frac{q-3}{2} \quad & \text{for odd $q$.}
\end{cases}
\end{equation*}
Next we define the set $\cB_q$ of \emph{forbidden blocks} as
\[
\cB_q:=
\begin{cases}
\{ \lbs \pm 1 \rbs\} \cup \bigcup_{m=1}^\infty \{ \lbs\pm 2, \pm m\rbs\}
  & \text{for $q=3$}, \\
\{ \lbs (\pm 1)^{h_q+1} \rbs\} \cup \bigcup_{m=1}^\infty \{ \lbs(\pm 1)^{h_q}, \pm m\rbs\}
  & \text{for even $q$ and} \\
\{ \lbs (\pm 1)^{h_q+1} \rbs\} \\
\; \cup \bigcup_{m=1}^\infty \{ \lbs(\pm 1)^{h_q},\pm 2,(\pm 1)^{h_q},\pm m\rbs\}
  & \text{for odd $q \geq 5$.}
\end{cases}
\]
The choice of the sign is the same within each block.
For example $\lbs 2,3 \rbs$, $\lbs -2,-3 \rbs \in \cB$ and $\lbs 2,-3 \rbs \not\in \cB$ for $q=3$.

We call a sequence $\lbs a_1,a_2,a_3,\ldots \rbs$ \emph{$q$-regular} if $\lbs a_k,a_{k+1}, \ldots, a_l \rbs \not\in \cB_q$ for all $1 \leq k < l$ and \emph{dual $q$-regular} if $\lbs a_l,a_{l-1}, \ldots, a_k \rbs \not\in \cB_q$ for all $1 \leq k < l$.
Denote by $\Ar$ and $\Adr$ the set of infinite $q$-regular and dual $q$-regular sequences $(a_i)_{i\in\IN}$, respectively.
\smallskip

A \emph{nearest $\lambda_q$-multiple continued fraction}, or \CF, is a formal expansion
\begin{equation}
\label{C1.2}
[a_0; a_1,a_2,a_3,\ldots]:= a_0 \lambda_q + \frac{-1}{a_1\lambda_q + \frac{-1}{a_2\lambda_q + \frac{-1}{a_3\lambda_q + \ldots } } }
\end{equation}
with $a_i \in \IZ\setminus \{0\}$, $i \geq 1$ and $a_0 \in \IZ$.

A \CF\, $[a_0; a_1,a_2,a_3,\ldots]$ is said to \emph{converge} if either $ [a_0; a_1,a_2,a_3,\ldots,a_l]$ has finite length or $\lim_{l \to \infty} [a_0; a_1,a_2,a_3,\ldots, a_l]$ exists in $\IR$. The notations for sequences, as introduced above, are also used for \CF's.


We say that a \CF \, is \emph{regular} or \emph{dual regular}, depending on whether the sequence $\lbs a_1,a_2,a_3,\ldots \rbs$ is $q$-regular or dual $q$-regular.
Regular and dual regular \CF\,'s are denoted by $\lb a_0; a_1,\ldots \rb$ and $\lbd a_0; a_1,\ldots \rbd$, respectively.

It follows from \cite[Lemmas 16 and 34]{MS08} that regular and dual regular \CF\,'s converge.
Moreover, it is known \cite{MS08} that $x$ has a regular expansion $x = \lb 0;a_1,a_2,\ldots \rb$ with leading $a_0=0$ if and only if $x \in \Iq:= \left[-\frac{\lambda_q}{2}, \frac{\lambda_q}{2} \right]$.


Convergent \CF\,'s can be rewritten in terms of the generators of the Hecke triangle group $G_q$:
if the expansion~\eqref{C1.2} is finite it can be written as follows
\begin{align*}
[a_0; a_1, a_2, a_3, \ldots, a_l]
&= a_0 \lambda_q + \frac{-1}{a_1\lambda_q + \frac{-1}{a_2\lambda_q + \frac{-1}{a_3\lambda_q + \ldots \frac{-1}{a_l\lambda_q} } }} \\
&= T^{a_0} \, ST^{a_1} \, ST^{a_2} \, ST^{a_3} \, \cdots \, ST^{a_l}\,0,
\end{align*}
since $\frac{-1}{a\lambda_q + x} = ST^a x$.
For infinite convergent \CF\,'s the expansion has to be interpreted as
\begin{align*}
[a_0; a_1, a_2, a_3, \ldots]
&= \lim_{l \to \infty}[a_0; a_1, a_2, a_3, \ldots, a_l] \\
&= \lim_{l\to \infty} T^{a_0} \, ST^{a_1} \, ST^{a_2} \, ST^{a_3} \, \cdots \, ST^{a_l}\, 0 \\
&= T^{a_0} \, ST^{a_1} \, ST^{a_2} \, ST^{a_3} \, \cdots \, 0.
\end{align*}

An immediate consequence of this is \cite[Lemma~2.2.2]{MM09}:
\begin{lemma}
\label{C1.1}
For a finite regular \CF\, one finds that
\begin{align*}
\lb a_0; a_1, \ldots,a_n,1^{h_q} \rb &= \lb a_0; a_1, \ldots,a_n-1,(-1)^{h_q} \rb
 \end{align*}
 for $q=2h_q+2$ and that
\begin{align*}
\lb a_0; \ldots,a_n,1^{h_q},2,1^{h_q} \rb &= \lb a_0; \ldots,a_n-1,(-1)^{h_q},-2,(-1)^{h_q} \rb
\end{align*}
 for $q=2h_q+3$.
\end{lemma}

\subsection{Special values and their expansions}
\label{C2}
The following results are well-known (see e.g.~\cite{MS08} and \cite[\S2.3]{MM09}).
The point $x=\mp\frac{\lambda_q}{2}$ has the regular \CF
\begin{equation}
\label{C2.2}
\mp\frac{\lambda_q}{2}
=
\begin{cases}
\lb 0;(\pm 1)^{h_q} \rb \qquad & \mbox{for even } q , \\
\lb 0;(\pm 1)^{h_q},\pm 2,(\pm 1)^{h_q} \rb \qquad & \mbox{for odd } q.
\end{cases}
\end{equation}

Define
\begin{align}
\label{C2.4}
R_q &:= r_q + \lambda_q
\quad \text{with} \\
\label{C2.5}
r_q &:=
\begin{cases}
\lb 0; \ov{1^{h_q-1},2} \rb & \text{for even $q$}, \\
\lb 0; \ov{3} \rb & \text{for $q = 3$,} \\
\lb 0; \ov{1^{h_q},2,1^{h_q-1},2} \rb & \text{for odd $q \geq 5$},
\end{cases}
\end{align}
whose expansion hence is periodic with period $\kappa_q$, where
\begin{equation}
\label{C2.3}
\kappa_q :=
\begin{cases}
h_q = \frac{q-2}{2}& \text{for even $q$} \\
2h_q+1 = q-2& \mbox{for odd $q$}.
\end{cases}
\end{equation}
The regular and dual regular $\lambda_q$-CF of the point $x=R_q$ has the form
\begin{align*}
R_q
&=
\begin{cases}
~ \lb 1 ; \ov{1^{h_q-1},2} \rb \quad & \mbox{for even $q$}, \\
~ \lb 1 ; \ov{3} \rb \quad & \mbox{for $ q=3$,}\\
~ \lb 1 ; \ov{1^{h_q},2,1^{h_q-1},2} \rb \quad & \mbox{for odd $q\geq 5$,  }
\end{cases} \\
\nonumber
&=
\begin{cases}
\lbd 0; (-1)^{h_q},\ov{-2,(-1)^{h_q-1}} \rbd & \text{for even $q$} , \\
\lbd 0; -2,\ov{-3} \rbd & \text{for $ q=3$,}\\
\lbd 0; (-1)^{h_q},\ov{-2,(-1)^{h_q},-2,(-1)^{h_q-1}} \rbd & \text{for odd $q\geq 5$. } \\
\end{cases}
\end{align*}
Moreover,
\begin{align}
&R_q = 1 \quad \text{and} \quad -R_q = S \, R_q
&& \text{for even $q$,} \label{R-relations-even}\\
&R_q^2 +(2-\lambda_q)R_q = 1 \quad \text{and} \quad -R_q = \big(T_{q}S\big)^{h_q+1} \, R_q
&& \text{for odd $q$},\label{R-relations-odd}
\end{align}
and $R_q$ satisfies the inequality
\[
\frac{\lambda_q}{2} < R_q \leq 1.
\]
\begin{remark}
\label{rem:r_and_minus_r}
Since $R_q=T_q \,r_q$ it follows from \eqref{R-relations-even} and \eqref{R-relations-odd} that on the one hand $-r_q=A_q r_q$ for some $A_q\in G_q$, but on the other hand it is clear
from \eqref{C2.5} that $r_q$ and $-r_q$ have different regular \CF~ expansions.

\end{remark}

\subsection{A lexicographic order}
\label{C3}
Let $x,y \in \Ir :=\left[ -R_q,R_q \right]$ have the regular \CF's $x = \lb a_0; a_1, \ldots \rb$ and $y = \lb b_0; b_1, \ldots \rb$.
Denote by $l(x)$ and $l(y)$ the number of entries (possibly infinite) in the above \CF's.
We introduce a \emph{lexicographic order} ``$\prec$'' for \CF's as follows:
If $a_i =b_i$ for all $0 \leq i \leq n$ and $l(x), l(y) \geq n$, we define
\[
x \prec y :\iff
\begin{cases}
a_0 < b_0 & \text{if $n=0$}, \\
a_n > 0 > b_n & \text{if $n>0$, both $l(x), l(y) \geq n+1$ and $a_n b_n < 0$}, \\
a_n < b_n & \text{if $n>0$, both $l(x), l(y) \geq n+1$ and $a_n b_n > 0$}, \\
b_n < 0 & \text{if $n>0$ and $l(x)=n$ }, \\
a_n > 0 & \text{if $n>0$ and $l(y)=n$}.
\end{cases}
\]
We also write $x \preceq y$ for $x \prec y$ or $x = y$.

This is indeed an order on regular \CF's, since \cite[Lemmas~22 and~23]{MS08} imply:
\begin{lemma}
\label{C3.1}
Let $x$ and $y$ have regular \CF's.
Then $x \prec y \iff x < y$.
\end{lemma}

\subsection{The generating interval maps $f_q$ and $f_q^\star$}
\label{C4}

The \emph{nearest $\lambda_q$-multiple map} $\nli{\cdot}$ is defined as
\[
\nli{\cdot} \colon \IR \to \IZ; \quad x \mapsto \nli{x}:= \nextinteger{ \frac{x}{\lambda_q} +
\frac{1}{2} }
\]
where $\nextinteger{\cdot}$ is the floor function
\[
\nextinteger{x} = n
\iff
\begin{cases}
n < x \leq n+1 \quad& \text{if } x > 0 , \\
n \leq x < n+1 \quad& \text{if } x \leq 0.
\end{cases}
\]
We also need the map $\nlid{\cdot}$ given by
\[
\nlid{\cdot} \colon \IR \to \IZ; \quad x \mapsto \nlid{x} :=
\begin{cases}
\nextinteger{ \frac{x}{\lambda_q} + 1 - \frac{R_q}{\lambda_q} } \qquad & \text{if } x \geq 0 , \\[5pt]
\nextinteger{ \frac{x}{\lambda_q} + \frac{R_q}{\lambda_q} } \qquad & \text{if } x < 0.
\end{cases}
\]
For the intervals
\[
\Iq = \left[ -\frac{\lambda_q}{2}, \frac{\lambda_q}{2} \right]
\quad \text{and} \quad
\Ir = \left[ -R_q, R_q \right]
\]
the interval maps $f_q \colon \Iq \to \Iq$ and $f_q^\star: I_{R_q} \to I_{R_q}$ are then defined as follows:
\begin{align}
\label{C4.2a}
f_q(x) &=
\begin{cases}
-\frac{1}{x} - \nli{\frac{-1}{x}} \lambda_q \qquad & \text{if } x \in \Iq \backslash \{0\},\\
0 & \text{if } x =0
\end{cases}
\intertext{and}
\label{C4.2b}
f_q^\star(y) &=
\begin{cases}
-\frac{1}{y} - \nlid{\frac{-1}{y}} \lambda_q \qquad & \text{if } y \in \Ir \backslash \{0\},\\
0 & \text{if } y =0.
\end{cases}
\end{align}
These maps generate the regular and dual regular \CF's in the following sense:

For given $x, y \in \IR$ the entries $a_i$ and $b_i$, $i \in \IZ_{\geq 0}$, in their \CF\, are determined by the algorithms:
\begin{itemize}
\item[(0)] $a_0 = \nli{x}$ and $x_1:= x-a_0\lambda_q \in \Iq$,
\item[(1)] $a_1 = \nli{\frac{-1}{x_1}}$ and $x_2:= \frac{-1}{x_1}-a_1\lambda_q = f_q(x_1) \in \Iq$,
\item[(\,$i$\,)] $a_i = \nli{\frac{-1}{x_i}}$ and $x_{i+1}:= \frac{-1}{x_i}-a_i\lambda_q = f_q(x_i) \in \Iq$ ,\, $i=2,3,\ldots$ ,
\item[($\star$)] the algorithm terminates if $x_{i+1}=0$,
\end{itemize}
and
\begin{itemize}
\item[(0)] $b_0 = \nlid{x}$ and $y_1:= y-b_0\lambda_q\in \Ir$,
\item[(1)] $b_1 = \nlid{\frac{-1}{y_1}}$ and $y_2:= \frac{-1}{y_1}-b_1\lambda_q= f_q^\star(y_1) \in \Ir$,
\item[(\,$i$\,)] $b_i = \nlid{\frac{-1}{y_i}}$ and $y_{i+1}:= \frac{-1}{y_i}-b_i\lambda_q= f_q^\star(y_i) \in \Ir$,\, $i=2,3,\ldots$ ,
\item[($\star$)] the algorithm terminates if $y_{i+1}=0$.
\end{itemize}
In \cite[Lemmas~17 and~33]{MM09} it is shown that these two algorithms lead to the regular and dual regular \CF's
\begin{equation*}
\label{C4.1}
x = \lb a_0; a_1, a_2, \ldots \rb
\quad \text{and} \quad
y = \lbd b_0; b_1, b_2, \ldots \rbd,
\end{equation*}
respectively.

\subsection{Markov partitions and transition matrices for $f_q$}
\label{C5}
It is known that the maps $f_q$ and $f_q^\star$ both possess the Markov property (cf. e.g.~\cite{MM09}).
This means that there exist partitions of the intervals $I_q$ and $I_{R_q}$ with the property that the set of boundary points is preserved by the map $f_q$ and $f_q^\star$, respectively.
Partitions with this property are called Markov partitions and we will now demonstrate how to construct these explicitly in this case (see also \cite[\S3.3]{MM09}).

Let $\mathcal{O}(x)$ denote the orbit of $x$ under the map $f_q$, i.e.
$$\mathcal{O}(x) = \left\{ f_q^n(x);\; n =0,1,2,\ldots \right\}.$$

We are interested in the orbits of the endpoints of $I_q$.
Because of symmetry it is enough to consider the orbit of $-\frac{\lambda_q}{2}$.
This orbit is finite; indeed
if $\# \{S\}$ denotes the cardinality of the set $S$ and $\kappa_q$ is given by \eqref{C2.3} then
$$
\# \{ \mathcal{O}\left(-\frac{\lambda_q}{2}\right)\} = \kappa_q +1.
$$
We denote the elements of $\mathcal{O}\left(-\frac{\lambda_q}{2}\right)$ by
\begin{align*}
\label{C5.1}
\phi_i &= f_q^i \left( - \frac{\lambda_q}{2} \right) = \lb 0;1^{h_q-i} \rb,
\qquad \, 0 \leq i \leq h_q=\kappa_q,\\
\intertext{for $q=2h_q+2 $ and by}
\phi_{2i} &= f_q^i \left( - \frac{\lambda_q}{2} \right) = \lb 0;1^{h_q-i},2, 1^{h_q}\rb, \quad 1\leq i \leq h_q ,\quad \text{and} \\
\nonumber
\phi_{2i+1} &= f_q^{h_q+i+1} \left( - \frac{\lambda_q}{2} \right) = \lb 0; 1^{h_q-i}\rb,
\quad 0 \leq i \leq h_q=\frac{{\kappa_q}-1}{2}
\end{align*}
for $q=2h_q+3$. Then the $\phi_i$'s can be ordered as follows (cf.~\cite{MS08})
\[
-\frac{\lambda_q}{2} = \phi_0 < \phi_1 < \phi_2 < \ldots < \phi_{{\kappa_q}-2} < \phi_{{\kappa_q}-1} = - \frac{1}{\lambda_q} < \phi_{{\kappa_q}} = 0.
\]
Define next $\phi_{-i}:=-\phi_i$, $0 \leq i \leq {\kappa_q}$.
The intervals
\begin{equation}
\label{Phi}
\Phi_i := \big[ \phi_{i-1},\phi_i \big]
\quad \text{and} \quad
\Phi_{-i}:= \big[ \phi_{-i} ,\phi_{-(i-1)} \big],
 \qquad 1 \leq i \leq {\kappa_q}
\end{equation}
define a Markov partition of the interval $\Iq$ for the map $f_{q}$. This means especially that
\begin{equation}
\label{MP}
\bigcup_{i\in A_{\kappa_q}} \Phi_{i} = \Iq
\quad \text{and} \quad
 \Phi_{i}^\circ \cap \Phi_{j}^\circ= \emptyset
\quad \text{for} \quad i \neq j \in A_{\kappa_q}
\end{equation}
holds. Here $A_{\kappa_q}=\left\{\pm1,\ldots,\pm \kappa_q\right\}$ and $S^\circ$ denotes the interior of the set $S$.
%

As in \cite{MM09} we introduce next a finer partition which is compatible with the intervals of monotonicity for $f_q$.

In the case $q=3$ where $\lambda_3=1$ define for $m=2,3,4,\ldots$ the intervals $J_m$ as
\begin{equation}
\label{C5.9}
J_2 = \left[ -\frac{1}{2}, -\frac{2}{5} \right]
\quad \text{and} \quad
J_m = \left[ -\frac{2}{2m-1}, -\frac{2}{2m+1} \right], \quad m = 3,4,\ldots
\end{equation}
and set $J_{-m}:= -J_m$ for $m=2,3,4,\ldots$.
This partition of $I_3$, which we denote by $\cM (f_3)$, is Markov.
The maps $f_3 |_{J_m}$ are monotone with $f_3 |_{J_m}(x) = -\frac{1}{x}-m$ and locally invertible with $(f_3 |_{J_m})^{-1}(y) = -\frac{1}{y+m}$ for $y \in f_3(J_m)$, $m=2,3,\ldots$.
\linebreak
For $q \geq 4$ consider the intervals $J_m$, $m=1,2,\ldots$, with
\begin{equation}
\label{C5.5}
\begin{split}
J_1 &= \left[ -\frac{\lambda_q}{2}, -\frac{2}{3\lambda_q} \right]
\quad \text{and} \\
J_m &= \left[ -\frac{2}{(2m-1)\lambda_q}, -\frac{2}{(2m+1)\lambda_q} \right], \quad m = 2,3,\ldots
\end{split}
\end{equation}
and set $J_{-m}:= -J_m$ for $m \in \IN$. The intervals $J_{m}$ are intervals of monotonicity for $f_q$, i.e the restriction $f_q |_{J_m}:x\mapsto -\frac{1}{x}-m \lambda_q$ is monotonically increasing.
Since some points in $\mathcal{O}\left(-\frac{\lambda_q}{2}\right)$ do not fall onto a boundary point of any of the intervals $J_m$, the partition given by these intervals has to be refined to become Markov.

For even $q$ define the intervals $J_{\pm1_i}$ as
\begin{equation}
\label{Jmq_eq_3}
J_{\pm 1_i} := J_{\pm 1} \cap \Phi_{\pm i},
\quad 1 \leq i \leq {\kappa_q},
\end{equation}
and therefore $J_{\pm 1_i}= \Phi_{\pm i}$ for $1 \leq i \leq {\kappa_q}-1$.
In this way one arrives at the partition $\cM (f_q)$, defined as
\begin{equation}
\label{Jmq_gt_3}
\Iq = \bigcup_{\epsilon = \pm}\left( \bigcup_{ i=1}^{{\kappa_q}} J_{\epsilon 1_i} \, \cup \,
\bigcup_{m=2}^\infty J_{\epsilon m}\right),
\end{equation}
which is clearly again Markov.

Consider next the case of odd $q \geq 5$. Here one has $\phi_{\pm i} \in J_{\pm 1}$ for $1 \leq i \leq {\kappa_q}-2$ and $\phi_{\pm ({\kappa_q}-1)} \in J_{\pm 2}$.
Hence define the intervals
\begin{align*}
\label{C5.8}
J_{\pm 1_i} &:= J_{\pm 1} \cap \Phi_{\pm i}
\quad 1 \leq i \leq {\kappa_q} -1\,
\text{and therefore}\, J_{\pm 1_i} =\Phi_{\pm i},\,
 1 \leq i \leq {\kappa_q} -2 ,\\
\nonumber
J_{\pm 2_ i} &:= J_{\pm 2} \cap \Phi_{\pm i}, \, i=\kappa_q-1,
\kappa_q.
\end{align*}
Then it is easy to see that the partition $\cM (f_q)$ defined by
\begin{equation*}
\label{C5.10}
\Iq = \bigcup_{\epsilon =\pm}\left(\bigcup_{ i=1}^{{\kappa_q}-1} J_{\epsilon 1_i} \,
\cup \bigcup_{i={\kappa_q}-1}^{{\kappa_q}} J_{\epsilon 2_i} \, \cup \, \bigcup_{ m=3}^\infty J_{\epsilon m} \right)
\end{equation*}
is a Markov partition.

A useful tool for understanding the dynamics of the map $f_q$ is the transition matrix $\IA = \big( \IA_{i,j} \big)_{i,j \in F_q}$, where $F_q$ is given by $\cM(f_q)$:
\[
F_q= \begin{cases}
\{ \pm 1_1, \ldots, \pm1_{\kappa_q-1},\pm 2, \pm 3,\ldots \} & \text{for even $q$,}\\
\{\pm 2, \pm 3,\ldots \} & \text{for $q=3$},\\
\{ \pm 1_1, \ldots, \pm1_{\kappa_q-1},\pm 2_{\kappa_q-1}, \pm 2_{\kappa_q}, \pm 3,\ldots \} & \text{for odd $q \geq 5$}.
\end{cases}
\]
Each entry $\IA_{i,j}$, $i,j \in F_q$ is given by
\[
\IA_{i,j} =
\begin{cases}
0 & \text{if } J_j^\circ \cap f_q(J_i^\circ) = \emptyset,\\
1 & \text{if } J_j^\circ \subset f_q(J_i^\circ).
\end{cases}
\]
As we will see in Chapter \ref{E}, the transition matrix is a crucial ingredient in our formula for the transfer operator.
The entries of the transition matrix are easy to obtain by keeping track of where the end points of the intervals $J_i$ are mapped to by $f_q$
and since $\IA_{i,j}=\IA_{-i,-j}$ for all $i,j\in F_q$, it is enough to consider the rows corresponding to ``positive'' indices.
A description of which entries of $\IA$ that are non-zero is given in the ensuing lemma (cf. also \cite{MM09}).
\begin{lemma}
For $q=3$ one finds that $\IA_{2,j}=1$ iff $j=-m$ for some $m\geq 2$ whereas $\IA_{i,j}=1$ for $i \geq 3$ and all $j\in F_q $.

For $q=2 h_q+2$ and $\kappa_q=h_q$ one finds that $\IA_{1_l,j}=1$ iff $j=1_{l+1}$ ($1\leq l \leq \kappa_q-1$).
Furthermore $\IA_{1_{\kappa_q-1},j}=1$ iff $j=m\geq 2$, whereas $\IA_{1_{\kappa_q},j}=1$ for all ``negative'' indices $j$ in $F_q$. Finally,
$\IA_{i,j}=1$ for all $i\geq 2$ and all $j\in F_q$.

For $q=2h_q+3 \geq 5$ and $\kappa_q=2h_q+1$ one finds that $\IA_{1_{2l-1},j}=1$ iff $j=1_{2l+1}$ ($1\leq l\leq h_q-1$).
Next, $\IA_{1_{2h_q-1},j}=1$ iff either $j=2_{\kappa_q}$ or $j=-m$ for some $m\geq 3$.
Furthermore $\IA_{1_{2l},j}=1$ iff $j=1_{2l+1}$ ($1\leq l\leq h_q-2$) and
$\IA_{1_{2h_q-2},j}=1$ iff $j=1_{2h_q}$ or $2_{\kappa_q}$. We also have $\IA_{1_{2h_q},j}=1$ for all ``negative'' indices $j$ in $F_q$.
Next one finds that $\IA_{2_{\kappa_q-1},j}=1$ iff $j=1_1$ and $\IA_{2_{\kappa_q},j}=1$ for all $j\in F_q, j\not= 1_1$.
Finally, $\IA_{i,j}=1$ for $i \geq 3$ and all $j\in F_q$.
\end{lemma}

Consider the local inverse $\vartheta_{n}: J_{n} \to \IR$ of the map $f_q$, defined on an interval of monotonicity $J_{n}$ given by (\ref{C5.9})
and (\ref{C5.5}) for $q=3$ and $q>3$, respectively.
This map can be expressed as
\begin{equation}
\label{LI}
\vartheta_{n}(x) := \left(f_q\big|_{J_{n}}\right)^{-1}(x) = \frac{-1}{x + n \lambda_q} = ST^{n} \, x.
\end{equation}
Its properties are given  in the following lemma.
\begin{lemma}
\label{C6.2}
The function $\vartheta_n$ has the following properties:
\begin{enumerate}
\item
it extends to a holomorphic function on $\IC \setminus \{-n\lambda_q\}$,
\item
it is strictly increasing on $(-\lambda_q, \lambda_q)$ and
\item
if either $0<n<m$ or $n<m<0$ or $m<0<n$ then $\vartheta_n (x) < \vartheta_m(x)$ for all $x \in (-\lambda_q, \lambda_q)$.
\end{enumerate}
\end{lemma}

\begin{proof}
The first property is straightforward and
since $\vartheta_n^\prime(x) = (n\lambda_q +x)^{-2}$ the derivative $\vartheta_n^\prime$ restricted to $(-\lambda_q, \lambda_q)$ is positive. Hence
 $\vartheta_n$ is strictly increasing on this interval.

To show the last property, consider the three cases separately and use that $-\lambda_q < x < \lambda_q$ to conclude that:
\begin{align*}
0 <n<m
&\iff 0 < n\lambda_q +x < m\lambda_q +x \,
\iff \frac{-1}{n\lambda_q +x} < \frac{-1}{m\lambda_q +x} \\
&\iff \vartheta_n(x) < \vartheta_m(x),
\end{align*}
\begin{align*}
n<m <0
&\iff n\lambda_q +x < m\lambda_q +x < 0\,
\iff \frac{-1}{n\lambda_q +x} < \frac{-1}{m\lambda_q +x} \\
&\iff \vartheta_n(x) < \vartheta_m(x) \qquad \text{and}
\end{align*}
\begin{align*}
m<0 < n
&\iff m\lambda_q +x < 0 < n\lambda_q +x \,
\iff \frac{-1}{m\lambda_q +x} > \frac{-1}{n\lambda_q +x} \\
&\iff \vartheta_m(x) > \vartheta_n(x).
\end{align*}
\end{proof}

\section{ The geodesic flow on Hecke surfaces}
\label{D}
\subsection{The unit tangent bundle}
The unit tangent bundle of $\IH$, which we denote by $\mathrm{UT} ( \IH)$, can be identified with $\IH \times \mathbb{S}^1$
where $\mathbb{S}^1$ is the unit circle.
A geodesic $\gamma$ on $\IH$ is either a half-circle based on $\IR$ or a line parallel to
the imaginary axis. Let $\Phi_t:\mathrm{UT} (\IH)\to \mathrm{UT}( \IH)$ be the geodesic flow along the oriented geodesic $\gamma$ and denote by $\phi_t$ the projection of $\Phi_t$ onto $\IH$.
The pair of base points of $\gamma$ are then denoted by $\gamma_\pm \in \IP_\IR$ where $\lim_{t \to \pm \infty} \phi_t = \gamma_{\pm}$.
An oriented geodesic $\gamma$ on $\IH$ is usually identified with the pair consisting of its base points $(\gamma_-,\gamma_+)$.

We often identify the unit tangent bundle of $\cM_q$, $\mathrm{UT}(\cM_q)$, with $\cF_q \times \mathbb{S}^1$.
Let $\pi_1^\star: \mathrm{UT}(\IH) \to \mathrm{UT}( \cM_q)$ be the extension of the projection map $\pi$ to $\mathrm{UT}(\IH)$.
Then the geodesic flow $\Phi_t$ on $\mathrm{UT}(\IH)$ projects to the geodesic flow $\pi_1^\star \circ \Phi_t$ on $\mathrm{UT}(\cM_q)$.
For simplicity we denote this by the same symbol, $\Phi_t$. The geodesic $\gamma^\star = \pi \gamma$ is a closed geodesic on $\cM_q$ if
and only if $\gamma_+$ and $\gamma_-$ are the two fixed points of a hyperbolic element in $G_q$.

\subsection{A Poincar\'e map for the geodesic flow and its associated Ruelle zeta function}
\label{D3}
 In \cite{MS08} a Poincar\'e section $\Sigma$ and a Poincar\'e map $\cP:\Sigma \to \Sigma$ for the geodesic flow on the Hecke surfaces were constructed.
To achieve this, the authors used properties of \CF\, expansions to construct a map $\tilde{\cP}:\tilde{\Sigma} \to \tilde{\Sigma}$ on a certain subset $\tilde{\Sigma} \subset \partial \cF_q \times \mathbb{S}^1$. The induced map $\cP:\Sigma \to \Sigma$ on the projection $\Sigma:= \pi_1^\star (\tilde{\Sigma} ) \subset \mathrm{UT}( \cM_q)$ was then shown to define a Poincar\'e map for the geodesic flow.

To be more precise, let $\gamma$ be a geodesic corresponding to an element $\tilde{z} \in \tilde{\Sigma}$ such that its base points $\gamma_\pm \in \IR$ have the regular and dual regular \CF \,expansions
\[
\gamma_- = \lb a_0; (\pm 1)^{k-1}, a_k, a_{k+1}, \ldots \rb
\quad \text{and} \quad
\gamma_+ = \lbd 0; b_1, b_2,\ldots \rbd.
\]
Then $\tilde{\cP}(\tilde{z})$ corresponds to the geodesic $g \gamma$ with base points $(g \, \gamma_-, g \, \gamma_+)$. Here $g \in G_q$ is determined by the property that the base points of the geodesic $g\gamma$ have the regular and dual regular expansions
\[
g\,\gamma_- = \lb a_k;a_{k+1}\ldots \rb
\quad \text{and} \quad
g\, \gamma_+ = \lbd 0; (\pm 1)^{k-1},a_0,b_1, b_2,\ldots \rbd
\]
corresponding to a $k-$fold shift of the bi-infinite sequence $$\ldots,b_2,b_1 \centerdot a_0,(\pm 1)^{k-1},a_k,a_{k+1},\ldots,$$ obtained by adjoining the dual regular and regular sequence corresponding to $\gamma_{+}$ and $\gamma_{-}$, respectively. It is possible to choose the sequence $b_1,b_2,\ldots$ such that the resulting bi-infinite sequence is \emph{regular}, i.e. that it does not contain \emph{any} forbidden blocks, not even across the ``zero marker'', $\centerdot$.
One can verify that the natural extension of the map $f_q$ is conjugated to the shift map on the space of regular bi-infinite sequences (cf. e.g.~\cite[Lemma 54]{MS08} or \cite{MM09}).
Indeed, if $a_0\neq 0$ then $\gamma_{-}\notin I_q$ but $S \, \gamma_- = -\gamma^{-1}_{-} \in \Iq$
and
\begin{equation}
\label{Poincare}
S\circ f_q^k \circ S \, \gamma_- = g\, \gamma_-.
\end{equation}
Our main interest lies in periodic orbits of the geodesic flow, that is, the flow along closed geodesics on $\mathcal{M}_q$.
Since $\cP$ is a Poincar{\'e} map, its periodic orbits correspond precisely to the periodic orbits of the geodesic flow.
From \eqref{Poincare} we see that there is a correspondence between the respective periodic orbits of the maps $\tilde{\cP}$ and $f_q$.
This correspondence is one-to-one and when projected to the surface, i.e. considering $\cP$ instead of $\tilde{\cP}$, it is bijective except for the periodic orbits under $f_q$ of the two points $\pm r_q$, which correspond to the same periodic orbit of $\cP$ (cf.~ Remark \ref{rem:r_and_minus_r}).
The periodic orbits of $f_q$ are determined by the points in $\Iq$ with periodic regular \CF~ expansions. Explicitly,
the base points of the closed geodesic $\gamma$, corresponding to the point $x = \lb 0; \ov{a_1, \ldots, a_n} \rb$ are given by
\[
\gamma_- = \lb a_1; \ov{a_2, \ldots, a_n,a_1} \rb
\quad \text{and} \quad
\gamma_+ = \lbd 0; \ov{a_n, \ldots, a_1} \rbd.
\]
Using (\ref{Poincare}) we conclude that a periodic orbit $\mathcal{O}^\star$ of $\cP$ can not have larger period than the corresponding periodic orbit of $f_q$, its period is in fact smaller if the \CF \, expansion of a point in the periodic orbit of $f_q$ contains a $1$ or $-1$.

A \emph{prime} (or \emph{primitive}) periodic orbit of the geodesic flow corresponds to a closed geodesic traversed once.
The analogous notion applies to closed geodesics, periodic orbits of $\cP$ and periodic orbits (and points) of $f_q$.
In particular, $x\in I_q$ is a prime periodic point of the map $f_q$ with period $n$ if $f_q^m(x)\not= x$ for all $0<m<n$ and $f_q^n(x)=x$.
Consider now a prime periodic orbit $\gamma^\star=(\gamma_-,\gamma_+)$ of the geodesic flow, determined by the prime periodic point $x^\star =
S \, \gamma_- = \lb 0 ; \ov{a_1, \ldots, a_n} \rb \in \Iq$ of the map $f_q$.
If the geodesic flow is appropriately normalized, the period $l(\gamma^\star)$ of $\gamma^\star$ is given by the length of the geodesic $\gamma$, which in turn is given by
\[
l(\gamma^\star) = 2 \ln \lambda
\]
where $\lambda$ is the larger one among the two real positive eigenvalues of the hyperbolic element
 $g^\star = ST^{a_1} ST^{a_2}\cdots ST^{a_n}\in G_q$, whose attracting fixed point is $x^\star \in \Iq$.
It is easy to verify that $f_q^n$ is in this case given precisely by $g^\star$ and a straightforward calculation shows that
\begin{align}
 \label{l-gamma-f}
l(\gamma^\star)
&= \ln \abs{\frac{\dd}{\dd x} f_q^n(x^\star)} = \sum_{l=0}^{n-1} \ln \abs{\frac{\dd}{\dd x}f_q\left( f_q^l(x^\star)\right)}
= \sum_{l=0}^{n-1} r \left( f_q^l(x^\star) \right),
\end{align}
where $r(x)=\ln f_q^\prime (x)$ and we used that $f_q^\prime(x)=\frac{1}{x^2}$ is positive and greater than one for $x$ in $I_q\setminus \{0\}$. Since $\tilde{\cP}^k( \tilde{z}^\star) = \tilde{z}^\star$ for some $k\leq n$ and
 $\tilde{z}^\star\in \tilde{\Sigma}$ corresponding to $x^\star$, the period $l(\gamma^\star)$ can also be written as
\begin{equation}
\label{l-gamma-P}
l(\gamma^\star) = \sum_{l=0}^{k-1} r\left( \tilde{\cP}^l(\tilde{z}^\star)\right)
\end{equation}
where $r(z^{\star}):=r(x^{\star})$.
Observe here that $r$ is not precisely the recurrence time function for the Poincar\'e map $\tilde{\cP}$ (cf.~\cite[Prop.~84]{MS08}).
However, when adding all pieces, the differences cancel and we nevertheless obtain the correct length.
To simplify the notation for later, if the periodic orbit $\mathcal{O}$ of $f_q$ and $\tilde{\mathcal{O}}$ of $\tilde{\cP}$ correspond to the closed geodesic $\gamma^\star$
we set $r_{\mathcal{O}}=r_{\tilde{\mathcal{O}}} = l(\gamma^\star)$. We remark here that $l(\gamma^\star)$ as given by \eqref{l-gamma-f} or \eqref{l-gamma-P} neither
depends on the choice of $x^\star \in \mathcal{O}$ in the first case nor of $\tilde{z}^\star \in \tilde{\mathcal{O}}$ in the second.
The Ruelle zeta function, $\zeta_R$, for the map $f_q$ is given by
\[
\zeta_R(s) = \exp \left( \sum_{n=1}^\infty \frac{1}{n} Z_n(s) \right)
\]
with
\begin{equation}
\label{Ruellez}
Z_n(s) = \sum_{x \in \Fix f_q^n} \exp \left(-s \sum_{k=0}^{n-1} \ln \, f_q^\prime \big( f_q^k (x) \big)\right).
\end{equation}
It is well-known that for $\re s$ large enough the above sums converge. Hence $\zeta_R(s)$
represents a holomorphic function in a half-plane of the form $\re s>\sigma$ for some $\sigma>0$.
A prime periodic orbit
\[
\mathcal{O} = \big( x, f_q(x), \ldots, f_q^{n-1}(x)\big)
\]
of period $n$ clearly contributes to all partition functions $Z_{ln}(s)$ with $l \in \IN$. Hence, if
$Z_{\mathcal{O}}(s)=\sum_{l=1}^\infty \frac{1}{l} \exp (-s l r_\mathcal{O})$ denotes the contribution of $\mathcal{O}$ to the Ruelle zeta function, one can use the Taylor expansion for $\ln (1-x)$ to see that
$$
\zeta_R(s)=\exp \left(\sum_\mathcal{O} Z_{\mathcal{O}}(s) \right)=\exp \left(-\sum_\mathcal{O} \ln \left(1-e^{-s\, r_\mathcal{O} }\right)\right),
$$
Therefore, summing over the set of all prime periodic orbits of $f_q$, leads to the well-known formula \cite{Ru94}
\begin{equation*}
\label{Ruelle zeta 1}
\zeta_R(s) = \prod_{\mathcal{O}}\left( 1 - e^{-sr_\mathcal{O}} \right)^{-1}.
\end{equation*}

Consider now the Ruelle zeta function for the map $\tilde{\cP}: \tilde{\Sigma} \to \tilde{\Sigma}$.
We know that the prime periodic orbits of this map and those of the map $f_q$ are in a one-to-one correspondence.
Furthermore, since $r_{\tilde{\mathcal{O}}}=r_{\mathcal{O}}$, all factors corresponding to the respective prime periodic orbits are equal.
We conclude that the Ruelle zeta function for $\tilde{\cP}$ is identical to that for $f_q$.

In \cite{MS08} it was shown that there is a one-to-one correspondence between the prime periodic orbits of the map
$\tilde{\cP}$ and the prime periodic orbits of the geodesic flow on $\mathrm{UT}( \cM_q)$, except for the two orbits
$\tilde{\mathcal{O}}_\pm$ determined by the endpoints $\big(S\, (\pm r_q), \mp r_q \big)$.
These two orbits coincide under the projection $\pi_q^\star: \mathrm{UT}( \IH) \to \mathrm{UT}( \cM_q)$.
However, the contributions of both of these two orbits are contained in the Ruelle zeta function $\zeta_R$.
The period of $\mathcal{O}_+$, the orbit of the point $r_q$ under the map $f_q$ is equal to $\kappa_q$ given by \eqref{C2.3}.
Define therefore partition functions $Z_n^{\mathcal{O}_+}(s)$, $n \in \IN$ as follows:
\begin{align*}
 Z_n^{\mathcal{O}_+}(s) &= 0 \qquad \text{for all $n$ with $\kappa_q \nmid n$} ,\\
Z_n^{\mathcal{O}_+}(s) &= \kappa_q \, \exp\left(-sl \, \ln \left(f_q^{\kappa_q} \right)^\prime (r_q)\right), \qquad n=\kappa_q l,\, l=1,2,\ldots .
\end{align*}
Then
\[
\exp \left(- \sum_{n=1}^\infty \frac{1}{n} \, Z_n^{\mathcal{O}_+}(s) \right)
=
\exp \left( - \sum_{l=1}^\infty \frac{1}{l} \, e^{-s l\, r_{\mathcal{O}_+}} \right)
=
1 - e^{-s \, r_{\mathcal{O}_+}}.
\]
Hence the Ruelle zeta function $\zeta_R^\cP (s)$ for the Poincar\'e map $\cP: \Sigma \to \Sigma$ of the geodesic flow $\Phi_t: \mathrm{UT}( \cM_q) \to \mathrm{UT}( \cM_q)$ has the form
\begin{equation*}
\label{Ruelle zeta 3}
\zeta_R^\cP (s)
=
\prod_{\mathcal{O} \neq \mathcal{O}_+} \left(1-e^{-s \, r_\mathcal{O}} \right)^{-1}.
\end{equation*}

\subsection{The Selberg zeta function}
\label{D2}
The \emph{Selberg zeta function}, $Z_S(s)$, for the Hecke triangle group $G_q$ is defined as
\begin{align}
\label{SZ-one-orbit}
Z_S(s) &=
\prod_{k=0}^\infty \prod_{\gamma^\star \text{ prime}}
  \left(1-e^{-(s+k)l(\gamma^\star)}\right) \nonumber \\
\intertext{where the inner product is taken over all prime periodic orbits $\gamma^\star$ of the geodesic flow on $\mathrm{UT}(\cM_q)$. It is now clear that we can write $Z_S(s)$ as}
Z_S (s) &=
\prod_{k=0}^\infty \prod_{\tilde{\mathcal{O}} \neq \tilde{\mathcal{O}}_+} \left(1-e^{-(s+k)r_\mathcal{O}}\right), \nonumber \\
\intertext{where the inner product is over all prime periodic orbits $\tilde{\mathcal{O}}$ of $\tilde{P}$, except for $\tilde{\mathcal{O}}_+ $.
For $\re{s}>1$ it can also be written in the form }
Z_S(s) &=
\frac{\prod_{k=0}^\infty \prod_{\mathcal{O}} \left(1-e^{-(s+k)r_\mathcal{O}}\right)}{\prod_{k=0}^\infty \left(1-e^{-(s+k)r_{\mathcal{O}_+}}\right)},
\end{align}
where the product $\prod_\mathcal{O}$ is over all prime periodic orbits $\mathcal{O}$ of the map $f_q$. If the period of $\mathcal{O}$ is equal to $l$ then
 $ r_{\mathcal{O}} = \ln \left( f_q^{l} \right)^\prime (x)$ for any $x\in\mathcal{O}$. In particular, $r_{\mathcal{O}_+} = \ln \left(f_q^{\kappa_q} \right)^\prime( r_q)$.
Note that the zeros of the denominator  of (\ref{SZ-one-orbit}) all lie in the left half $s$-plane. We will show next how $Z_S(s)$ can be expressed in terms of the Fredholm determinant of the
 transfer operator for the map $f_q:I_q\to I_q$.

\section{Ruelle's transfer operator for the map $f_q$}
\label{E}
If $g:I_q\to \mathbb{C}$ is a function on the interval $I_q$ then Ruelle's transfer operator for the map $f_q$, $\mathcal{L}_s$, acts on $g$ as follows:
\begin{equation}
\label{E1.2}
\mathcal{L}_s g(x) = \sum_{y\in f_q^{-1}(x)} e^{-s \,r(y)} \, g(y)
\end{equation}
where $r(y)=\ln f_q^\prime (y)$ and $\re{s} >1$ to ensure convergence of the series.
To get a more explicit form for $\mathcal{L}_s$ one has to determine the set of preimages $f_q^{-1}(x)$ of an arbitrary point $x
\in \Iq$.
For this purpose recall the Markov partition $\Iq = \bigcup_{i\in A_{\kappa_q}}\Phi_i$ from (\ref{MP}) as well as the local inverses of $f_q$, $\vartheta_n$ from (\ref{LI}).
Using the Markov property of the map $f_q$, the preimages of points in $\Phi_i^\circ$ can be characterized by the following lemma:
\begin{lemma}
For $x\in\Phi_i^\circ \subset I_q $ and $i\in A_{\kappa_q}$ the preimage $f_q^{-1}(x)$ is given by the
set $f_q^{-1}(x)=\{y\in I_q:y=\vartheta_n(x),\, n\in\cN_i \}$ with $\cN_i=\bigcup _{j\in
A_{\kappa_q}}\cN_{i,j}$ and $\cN_{i,j}=\{n\in \mathbb{Z}: \vartheta_n(\Phi_i^\circ)\subset\Phi_j^\circ\}$
\end{lemma}
\begin{proof}
The preimages of a point $x$ in the open interval $\Phi_i^{\circ}$ can be determined from its
\CF\,expansion $x=\lb 0;a_1,a_2,\ldots\rb$. The boundary points of the interval $\Phi_i$ belong to the orbit $\mathcal{O}\left(-\frac{\lambda_q}{2}\right)$.
To be precise, $\Phi_i=[\lb 0;1^{h_q+1-i}\rb,\lb 0;1^{h_q-i}\rb ],\, 1\leq i\leq h_q$ for $q=2h_q+2$ and
$\Phi_{2i+1}=[\lb 0;1^{h_q-i},2,1^{h_q}\rb,\lb 0;1^{h_q-i}\rb ],\, 1\leq i \leq h_q$, and
$\Phi_{2i}=[\lb 0;1^{h_q+1-i}\rb,\lb 0;1^{h_q-i},2,1^{h_q}\rb ],$ $1\leq i \leq h_q$ for $q=2h_q+3$. For the intervals
$\Phi_{-i}$ one gets analogous expressions with negative entries in the \CF\, expansions. Using the lexicographic order  one concludes that the \CF\, expansion of a point $x\in
\Phi_i^\circ$ for $q=2h_q+2$ must be either of the form $x=\lb
0;1^{h_q+1-i},-m,\ldots\rb$ with $m\geq 1$ or of the form $x=\lb 0;1^{h_q-i},m,\ldots\rb$ with
$m\geq 2$. It is easy to see that the set of $n\in \mathbb{Z}$ such that $\vartheta_n(x)=\lb 0;n,x\rb \in I_q$ only depends on the interval $\Phi_i$ with $x\in \Phi_i$. Here $\lb 0;n,x\rb$ denotes the concatenation of the corresponding sequences. Furthermore, if $\vartheta_n(x)\in \Phi_j^\circ$ for some $x\in \Phi_i^\circ$ then $\vartheta_n(\Phi_i^\circ)\subset \Phi_j^\circ$ and hence $\cN_i=\bigcup _{j\in
A_{\kappa_q}}\cN_{i,j}$. The same reasoning applies to the case $q=2h_q+3$.
\end{proof}
\begin{remark}
We will define the operator $\mathcal{L}_s$ on a space of piecewise continuous functions, hence
it is enough to determine the preimages of points in the interior $\Phi_i^\circ$ of the intervals $\Phi_i$.
In general, points on the boundary of an interval $\Phi_i$ can have more preimages than points in the
interior.
\end{remark}
For $n=1,2,\ldots$ define $\mathbb{Z}_{\geq n}:=\{l\in \mathbb{N}: l\geq n\}$ and $\mathbb{Z}_{\leq -n}:=\{l\in \mathbb{Z}: l\leq -n\}$.
We are now able to determine the sets $\cN_{i,j}$ explicitly.
\begin{lemma}
The sets $\cN_{i,j}=\{n\in \mathbb{Z}: \vartheta_n(\Phi_i)\subset\Phi_j\}$ are given by the following expressions.
For $q=2h_q+2$ we have 
\begin{align*}
\cN_{1,h_q}&=\mathbb{Z}_{\geq 2},\, \cN_{1,-h_q}=\mathbb{Z}_{\leq -1},\nonumber\\
\cN_{i,i-1}&=\{1\}, \, \cN_{i,h_q}=\mathbb{Z}_{\geq 2}, \, \cN_{i,-h_q}=\mathbb{Z}_{\leq -1},\, 2\leq i\leq h_q.
\end{align*}
For $q=3$ we have: 
\begin{align*}
\cN_{1,1}&=\mathbb{Z}_{\geq 3},\, \cN_{1,-1}=\mathbb{Z}_{\leq -2}.
\end{align*}
For $q=2h_q+3\ge5$ we have:
\begin{align*}
\cN_{1,2h_q}&=\{2\}, \,\cN_{1,-2h_q}=\{-1\},\, \cN_{1,2h_q+1}= \mathbb{Z}_{\geq 3}, \,
\cN_{1,-(2h_q+1)}=\mathbb{Z}_{\leq -2},\nonumber\\
\cN_{2,-2h_q}&=\{-1\}, \,\cN_{2,2h_q+1}=\mathbb{Z}_{\geq 2}, \,\cN_{2,-(2h_q+1)}=\mathbb{Z}_{\leq -2},\\
\cN_{i,i-2}&=\{1\}, \, 3\leq i\leq \kappa_q,\, \cN_{i,-2h_q}=\{-1\}, \, 1\leq i \leq \kappa_q,\nonumber\\
\cN_{i,2h_q+1}&=\mathbb{Z}_{\geq 2}, \,\cN_{i,-(2h_q+1)}=\mathbb{Z}_{\leq -2},\, 1\leq i \leq \kappa_q \nonumber.
\end{align*}
Furthermore, $\cN_{-i,j}=-\cN_{i,-j}$ for $ i,j\in A_{\kappa_q}$ and
for all pairs of indices not listed above, the set $\cN_{i,j}$ is empty.
\end{lemma}
\begin{proof}
We will give a proof for the case $q=2h_q + 2$. The proof for the case of odd $q$ is similar.
If $x\in \Phi_1^\circ$ then
either $x=\lb 0;1^{h_q},-m,\ldots\rb$ for some $m\geq 1$ or $x=\lb 0;1^{h_q-1},m,\ldots\rb$ for some
$m\geq 2$. In both cases $\lb 0;1,x\rb\notin I_q$ whereas $\lb 0;-1,x\rb\in \Phi_{-h_q}$
and $\lb 0;\pm n,x\rb\in \Phi_{\pm h_q}$ for $n\geq 2$.
For $x\in \Phi_{h_q}^\circ$ one has
$x=\lb 0;1,-m,\ldots\rb$ for some $m\geq 1$ or $x=\lb 0;m,\ldots\rb$ for some $m\geq 2$.
In this case $\lb 0;1,x\rb\in \Phi_{h_q-1}$, $\lb 0;\pm n,x\rb\in \Phi_{\pm h_q}$ and $\lb 0;-1,x\rb\in
\Phi_{-h_q}$. Finally, for $x\in \Phi_i^\circ,\, 2\leq i\leq h_q-1=\kappa_q-1$ one has either $x=\lb
0;1^{h_q+1-i},-m,\ldots\rb$ for some $m\geq 1$ or $x=\lb 0;1^{h_q-i},m,\ldots\rb$ for some $m\geq 2$. In
both cases $\lb 0;1,x\rb\in \Phi_{i-1}$, $\lb 0;\pm n,x\rb\in \Phi_{\pm h_q}$ for all $n\geq 2$ and
$\lb 0;-1,x\rb\in \Phi_{h_q}$.
It is clear that the set $\cN_{i,j}$ is empty for all combinations of indices which are not listed above.
That $\cN_{-i,j}=-\cN_{i,-j}$ for all $i,j\in A_{\kappa_q}$ follows immediately from the fact that $\Phi_{-i}=-\Phi_i$ in combination with the explicit form of the map $\vartheta_n$.
\end{proof}
The previous lemma allows us to derive explicit expressions for the transfer operator.
Using the index sets $\cN_{i}=\bigcup_{j\in A_{\kappa_q}}\cN_{i,j}$ we can rewrite $\mathcal{L}_s$ in \eqref{E1.2} as
\begin{equation*}
\label{E1.3}
\mathcal{L}_s g(x)= \sum_{i\in A_{\kappa_q}} \chi_{\Phi_i}(x) \sum_{n\in\cN_i}
\big(\vartheta_n^\prime(x)\big)^s \, g\big(\vartheta_n(x)\big),
\end{equation*}
where $\chi_{\Phi_i}$ is the characteristic function of the set $\Phi_i$.
If we now introduce vector valued functions $\underline{g}=(g_i)_{i\in A_{\kappa_q}}$
with $g_i:= g|_{ \Phi_i}$, then the operator $\mathcal{L}_s$ can also be written as follows
\begin{align*}
(\mathcal{L}_s \underline{g})_i(x)
&= \sum_{j\in A_{\kappa_q} }\sum_{n\in\cN_{i,j}} \big(\vartheta_n^\prime(x)\big)^s \,
g_j\big(\vartheta_n(x)\big) \\
&= \sum_{j\in A_{\kappa_q} } \sum_{n\in\cN_{i,j}}
\left(\frac{1}{z+n\lambda_q}\right)^{2s} \, g_j \left(\frac{-1}{z+n\lambda_q} \right), \quad x \in
\Phi_i.
\end{align*}
If $g_i$ is continuous on $\Phi_i$ for all $i\in A_{\kappa_q}$ then $(\mathcal{L}_s \underline{g})_i$ is also continuous on
$\Phi_i$ since $\vartheta_n(\Phi_i^\circ)\subset \Phi_j^\circ$ for $n\in\cN_{i,j}$. This implies that $\mathcal{L}_s$ is well defined on the Banach space,
$B=\oplus_{i\in A_{\kappa_q}} C(\Phi_i)$, of piecewise continuous functions on the intervals $\Phi_i$.
To be able to give explicit expressions for $\mathcal{L}_s$ on the space $B$ we need two auxiliary operators $\mathcal{L}_{\pm n,s}^\infty$ and $\mathcal{L}_{\pm n,s}$. For $n\in \mathbb{N}$ these are defined by
\begin{equation}
\label{aux}
\mathcal{L}_{\pm n,s}^\infty g(x) =\sum_{l=n}^\infty \frac{1}{(x\pm l\lambda_q)^{2s}}\,g\left(\frac{-1}{x\pm l \lambda_q}\right),\\
\end{equation}
and by
\begin{equation}\label{aux1}
\mathcal{L}_{\pm n,s}g(x) =\frac{1}{(x\pm n\lambda_q)^{2s}}\,g\left(\frac{-1}{x\pm n \lambda_q}\right).
\end{equation}
Then we have
\begin{lemma}
For $q=3$ the operator $\mathcal{L}_s$ is given by
\begin{align*}
(\mathcal{L}_s\underline{g})_1 &= \mathcal{L}_{3,s}^\infty g_1 +\mathcal{L}_{-2,s}^\infty g_{-1},\nonumber\\
(\mathcal{L}_s\underline{g})_{-1} &= \mathcal{L}_{2,s}^\infty g_1 +\mathcal{L}_{-3,s}^\infty g_{-1}.
\end{align*}
For $q=2h_q+2$ one has
\begin{align*}
(\mathcal{L}_s\underline{g})_1 &=\mathcal{L}_{2,s}^\infty g_{h_q}+\mathcal{L}_{-1,s}^\infty g_{-h_q},\nonumber\\
(\mathcal{L}_s\underline{g})_i &= \mathcal{L}_{1,s} g_{i-1}+ \mathcal{L}_{2,s}^\infty g_{h_q}+ \mathcal{L}_{-1,s}^\infty g_{-h_q},\, 2\leq i\leq h_q,
\end{align*}
and
\begin{align*}
(\mathcal{L}_s\underline{g})_{-1} &=\mathcal{L}_{1,s}^\infty \, g_{h_q}+\mathcal{L}_{-2,s}^\infty \, g_{-h_q},\nonumber\\
(\mathcal{L}_s\underline{g})_{-i} &= \mathcal{L}_{-1,s} \, g_{-(i-1)}+ \mathcal{L}_{1,s}^\infty \, g_{h_q}+ \mathcal{L}_{-2,s}^\infty \, g_{-h_q},\, 2\leq i\leq h_q.
\end{align*}
For $q=2h_q+3$ one has
\begin{align*}
(\mathcal{L}_s\underline{g})_{1} &=\mathcal{L}_{2,s} \, g_{2h_q}+\mathcal{L}_{3,s}^\infty \, g_{2h_q+1}+\mathcal{L}_{-2,s}^\infty \, g_{-(2h_q+1)}+\mathcal{L}_{-1,s} \, g_{-2h_q},\nonumber\\
(\mathcal{L}_s\underline{g})_{2}&=\mathcal{L}_{2,s}^\infty \, g_{2h_q+1}+ \mathcal{L}_{-2,s}^\infty \, g_{-(2h_q+1)} + \mathcal{L}_{-1,s} \, g_{-2h_q},\\
(\mathcal{L}_s\underline{g})_{i} &= \mathcal{L}_{1,s} \, g_{i-2}+ \mathcal{L}_{2,s}^\infty \, g_{2h_q+1}+\mathcal{L}_{-2,s}^\infty \, g_{-(2h_q+1)}+ \mathcal{L}_{-1,s} \, g_{-2h_q},\, 1\leq i\leq 2h_q+1,\nonumber
\end{align*}
and
\begin{align*}
(\mathcal{L}_s\underline{g})_{-1} &=\mathcal{L}_{1,s} \, g_{2h_q}+\mathcal{L}_{2,s}^\infty \, g_{2h_q+1}+\mathcal{L}_{-3,s}^\infty \, g_{-(2h_q+1)}+\mathcal{L}_{-2,s} \, g_{-2h_q},\nonumber\\
(\mathcal{L}_s\underline{g})_{-2}&= \mathcal{L}_{1,s} \, g_{2h_q}+ \mathcal{L}_{2,s}^\infty \, g_{2h_q+1} + \mathcal{L}_{-2,s}^\infty \, g_{-(2h_q+1)},\\
(\mathcal{L}_s\underline{g})_{-i} &= \mathcal{L}_{1,s} \, g_{2h_q}+ \mathcal{L}_{2,s}^\infty \, g_{2h_q+1}+\mathcal{L}_{-2,s}^\infty \, g_{-(2h_q+1)}+\mathcal{L}_{-1,s} \, g_{2-i},\, 1\leq i\leq 2h_q+1.\nonumber
\end{align*}
\end{lemma}

Unfortunately the operator $\mathcal{L}_s$ is not of trace class on the space of piecewise continuous functions. In fact, it is even not compact on this space.

Much better spectral properties however can be achieved by defining $\mathcal{L}_s$ on the Banach space
$B=\oplus_{i\in A_{\kappa_q}} B(D_i)$ with $B(D_i)$ the Banach space of holomorphic functions on a
certain open disc $D_i\subset \mathbb{C}$ with $\Phi_i \subset D_i$, for all $i \in A_{\kappa_q}$, and continuous on the closed disc $\overline{D}_i$,
 together with the sup norm.
This is possible since all the maps $\vartheta_{\pm m}, \, m\geq 1$ have holomorphic extensions to the disks, whose existence is given by the following lemma.
\begin{lemma}
\label{E1.5}
There exist open discs $D_i\subset \IC$, $i\in A_{\kappa_q} $, with $\Phi_i\subset D_i$ and $\vartheta_n(\overline{D}_i)\subset D_j$ for all $n\in \cN_{i,j}$ .
\end{lemma}
For the proof of the Lemma it suffices to show the existence of open intervals $ I_i\subset \mathbb{\IR}$, $i \in A_{\kappa_q}$ with
\begin{itemize}
\item $\Phi_i \subset I_i$ and
\item $\vartheta_n( \ov{I_i}) \subset I_j$ for all $n \in \cN_{i,j}$.
\end{itemize}

Since the maps $\vartheta_n$ are conformal it is clear that the discs $D_i$ with center on the real axis and intersection equal to the open intervals $I_i$ then satisfy Lemma \ref{E1.5}.

Using (\ref{LI}) the two conditions on $I_i$ can also be written as
\begin{equation}
\label{E1.9}
\Phi_i \subset I_i \quad \text{and} \quad ST^n \, \overline{I}_i \subset I_j \quad \text{for all } n \in \cN_{i,j}\quad\text{and all } i,j\in A_{\kappa_q}.
\end{equation}

In the cases $q=3$ and $q=4$ we give explicit intervals fulfilling
conditions~\eqref{E1.9}.
For the case $q \geq 5$ we first show the existence of intervals $I_i$ satisfying \eqref{E1.9} with the second condition replaced
by the (slightly) weaker
\begin{equation*}
\label{E1.7}
ST^n \, I_i \subset I_j \quad \text{for all }n \in \cN_{i,j}.
\end{equation*}
The existence of intervals $I_i$ satisfying~\eqref{E1.9} then follows by a simple perturbation argument.

\begin{lemma}
\label{E3.3}
The intervals
\begin{align*}
\label{E3.3a}
&I_1:=\left(-1,\frac{1}{2}\right)\quad \text{and}\quad
I_{-1}:=-I_1 \quad\text{for}\quad q=3,
\end{align*}
and
\begin{align*}
&I_{1}:= \left( -1,\frac{\lambda_q}{4} \right)\quad \text{and} \quad I_{-1}:=-I_1 \quad \text{for}\quad q=4
\end{align*}
satisfy the conditions~\eqref{E1.9} and hence Lemma \ref{E1.5} holds for $q=3$ and $4$.
\end{lemma}

\begin{proof}
Since $\lambda_3 = 1$ and $\lambda_4 = \sqrt{2}$, the above intervals $I_i, i=\pm 1$ obviously satisfy
\begin{align*}
\Phi_1 = \left[-\frac{\lambda_q}{2},0\right] \subset I_1 \qquad \text{and}
\qquad \Phi_{-1} = \left[0, \frac{\lambda_q}{2}\right] \subset I_{-1}.
\end{align*}
For $q=3$ one has $\mathcal{N}_{1,1}=\mathbb{Z}_{\geq 3}$, $\mathcal{N}_{1,-1}=\mathbb{Z}_{\leq
-2}$, $\mathcal{N}_{-1,-1}=-\mathbb{Z}_{\geq 3}$ and $\mathcal{N}_{-1,1}=-\mathbb{Z}_{\leq
-2}$. Hence we have to show that $\theta_n(\overline{I_1})\subset I_1$ for all $n\geq 3$ and
$\theta_n(\overline{I_1})\subset I_{-1}$ for all $n\leq -2$. Since all maps involved are strictly
increasing, it is enough to show that $\theta_n(-1) > -1$ and $\theta_n(\frac{1}{2})< \frac{1}{2}$ for all
$n\geq 3$, as well as that $\theta_n(-1) > -\frac{1}{2}$ and $\theta_n(\frac{1}{2})< 1$ for all $n\leq -2$.
But this is not hard to show:
$\theta_n(-1) = \frac{-1}{-1+n}\geq \frac{-1}{2}> -1$ and
$\theta_{-n}(\frac{1}{2})=\frac{-1}{\frac{1}{2}-n}=\frac{1}{n-\frac{1}{2}}\leq \frac{2}{3}<1$ for all
$n\geq 2$. Furthermore, $\theta_n(-1)>0$ if $n<0$ and $\theta_n(\frac{1}{2}<0$ if $n>0$.
Since $\mathcal{N}_{-i,j}=-\mathcal{N}_{i,-j}$ and $I_{-1}=-I_1$ the result for the interval
$I_{-1}$ follows directly.

Consider now the case $q=4$. Then one has $\mathcal{N}_{1,1}=\mathbb{Z}_{\geq 2}$,
 $\mathcal{N}_{1,-1}=\mathbb{Z}_{\leq -1}$, $\mathcal{N}_{-1,-1}=-\mathbb{Z}_{\geq 2}$
and $\mathcal{N}_{-1,1}=-\mathbb{Z}_{\leq -1}$. One sees that $\theta_n(-1) =
\frac{-1}{-1+n\lambda_4}\geq
\frac{-1}{-1+2\lambda_4}> -1$ since $2 \lambda_4=2\sqrt{2}> 2$ and
$\theta_n(\frac{\lambda_4}{4})=\frac{-1}{\frac{\lambda_4}{4}+n\lambda_4}< 0<\frac{\lambda_4}{4}$ for all
$n\geq 2$. Furthermore $\theta_{-n}(-1) = \frac{-1}{-1-n\lambda_4}\geq 0>-\frac{\lambda_4}{4}$ and
$\theta_{-n}(\frac{\lambda_4}{4}) = \frac{-1}{\frac{\lambda_4}{4}-n\lambda_4}\leq
\frac{1}{\lambda_4-\frac{\lambda_4}{4}}=\frac{4}{3\lambda_4}<1$ for all $n\geq 1$ since $3 \sqrt2 > 4$.
Since $\mathcal{N}_{-i,j}=-\mathcal{N}_{i,-j}$ and $I_i=-I_{-i}$ the lemma is proved.
\end{proof}


To prove Lemma~\ref{E1.5} for $q \geq 5$ we will need the next four lemmas.
\label{E4}

\begin{lemma}
\label{E4.5}
For $q =2h_q+2, h_q\geq 2$ and $0\leq i\leq h_q$
\begin{equation*}
\label{E4.5a}
(ST)^{h_q-i} \, \left(-\frac{\lambda_q}{2}\right) = \lb -1;(-1)^i \rb, \qquad i =0,\ldots,h_q
\end{equation*}
and
\begin{equation}
\label{E4.5b}
\begin{split}
-\lambda_q
&= \lb -1; \rb < \lb -1;(-1)^1 \rb < \ldots \\
&< \lb -1;(-1)^{h_q-1} \rb < \lb -1;(-1)^{h_q} \rb = -\frac{\lambda_q}{2}.
\end{split}
\end{equation}
\end{lemma}

\begin{proof}
Using the $\lambda_q$-CF of $\lambda_q$ in ~\eqref{C2.2} we have
\begin{align*}
(ST)^{h_q-i} \, \left(-\frac{\lambda}{2}\right)
&= (ST)^{h_q-i} \, (ST)^{h_q} \, 0
 = (ST)^{-i-2} \, 0 \\
&= (T^{-1}S)^i \, T^{-1}\,ST^{-1}S \, 0
 = T^{-1}\, (ST^{-1})^i \,0 \\
&= \lb -1;(-1)^i \rb.
\end{align*}
This shows the first part of Lemma \ref{E4.5}.
By definition it is clear that $- \lambda_q=\lb -1; \rb$ and then inequalities \eqref{E4.5b} follow immediately from the lexicographic order in Section \ref{C3}.
\end{proof}
\begin{lemma}
\label{E4.4}
For $q=2h_q+2, h_q\geq 2 $ define the intervals
$I_i:=\left( \lb -1;(-1)^i \rb , \frac{\lambda_q}{4} \right)$ for $1\leq i \leq h_q$ and let $I_{-i}:=-I_i$. Then
\begin{align*}
&\vartheta_{\pm n} ( \overline I_{\pm i} ) \subset I_{\pm h_q} \quad\text{for all } n \geq 2,\, i=1,\ldots, h_q,\\
&\vartheta_{\pm n} ( \overline I_{\mp i} ) \subset I_{\pm h_q}\quad \text{for all } n \geq 1,\, i=1,\ldots,h_q , \\
&\vartheta_{\pm 1} ( I_{\pm i}) \subset I_{\pm {i-1}} \quad \text{for all } i=2,\ldots,h_q.
\end{align*}
Hence $\vartheta_n(I_i)\subset I_j$ for all $n \in \mathcal{N}_{i,j}$.
\end{lemma}

\begin{proof}
Since $\mathcal{N}_{i,h_q}=\mathbb{Z}_{\geq 2}$ for all $ 1\leq i\leq h_q$ we have to show that
$\vartheta_n ( \overline I_{ i} ) \subset I_{ h_q}$ for all $ 1\leq i\leq h_q$ and all $n\geq 2$. On the one hand, by Section \ref{C3} we see immediately that
$\vartheta_n (\lb -1;(-1)^i \rb)=\lb 0;n-1,(-1)^i \rb > -\frac{\lambda_q}{2}$ and on the other hand, $\vartheta_n(\frac{\lambda_q}{4})=\frac{-1}{n \lambda_q +
\frac{\lambda_q}{4}}<0<\frac{\lambda_q}{4}$. Hence $\vartheta_n(\overline{I}_i)\subset I_{h_q}$.
Consider next the case $\mathcal{N}_{i,-h_q}=\mathbb{Z}_{\leq -1}$ for $1\leq i\leq h_q$. There one has
$\vartheta_{-n}(\lb -1;(-1)^i \rb)=\lb 0;-n-1,(-1)^i \rb>0>-\frac{\lambda_q}{4}$.
Furthermore, for $n\ge 1$, one finds that $\vartheta_{-n}(\frac{\lambda_q}{4})=\frac{-1}{-n \lambda_q
+\frac{\lambda_q}{4}} \leq \frac{4}{3 \lambda_q}< \frac{\lambda_q}{2}$ since $\lambda_q\geq \sqrt{3}$ for $q\geq 6$.
Hence $\vartheta_{-n}(\overline{I}_i)\subset I_{-h_q}$ for all $n\geq 1$.
Consider finally the case $\mathcal{N}_{i,i-1}=\{1\}$ for $2\leq i\leq h_q$.
In this case $\vartheta_1(\lb -1;(-1)^i \rb)=\lb -1;(-1)^{i-1} \rb$ and since $\vartheta_1(\frac{\lambda_q}{4})=\frac{-1}{ \lambda_q+\frac{\lambda_q}{4}}<0<\frac{\lambda_q}{4}$
it follows that $\vartheta_1(I_i)\subset I_{i-1}$ for all $2\leq i\leq h_q$. The intervals $I_{-i}$ again have analogous properties.
\end{proof}

\begin{lemma}
\label{E4.2}
For $q=2h_q+3, h_q\geq 1$ one has
\begin{equation}
\label{E4.2a}
\begin{split}
(ST)^{h_q-i} \, \left(-\frac{\lambda_q}{2}\right) &= \lb -1;(-1)^i,-2,(-1)^{h_q} \rb \qquad \text{for} \quad 0\leq i\leq h_q, \\
(ST)^{h_q+1-i} \, \big(\lb 0;1^{h_q} \rb \big) &= \lb -1;(-1)^{i} \rb \qquad \text{for} \quad 1 \leq i\leq h_q
\end{split}
\end{equation}
and
\begin{equation}
\label{E4.2b}
\begin{split}
-\lambda_q & =
 \lb -1; \rb
< \lb -1;-2,(-1)^{h_q} \rb < \lb -1;-1 \rb < \lb -1;(-1)^1,-2,(-1)^{h_q} \rb \\
& < \lb -1;(-1)^2 \rb
< \lb -1;(-1)^2,-2,(-1)^{h_q} \rb < \lb -1;(-1)^3 \rb < \ldots\\
& < \lb -1;(-1)^{h_q-1},-2,(-1)^{h_q} \rb < \lb -1;(-1)^{h_q} \rb < \lb -1;(-1)^{h_q},-2,(-1)^{h_q} \rb \\
& = -\frac{\lambda_q}{2}.
\end{split}
\end{equation}
\end{lemma}

\begin{proof}
Using the $\lambda_q$-CF in ~\eqref{C2.2} it is easy to see that
\begin{align*}
(ST)^{h_q-i} \, \left(-\frac{\lambda}{2}\right)
&= (ST)^{h_q-i} \, (ST)^{h_q} \, ST^2 \, (ST)^{h_q}\, 0 \\
&= (ST)^{-i-3} \, ST^2 \, (ST)^{h_q}\, 0 \\
&= T^{-1} \, (ST^{-1})^i \,ST^{-1}\,ST^{-1}S \, ST^2 \, (ST)^{h_q}\, 0 \\
&= T^{-1} \, (ST^{-1})^i \,ST^{-1}\,(ST)^{h_q+1} \, 0 \\
&= T^{-1} \, (ST^{-1})^i \,ST^{-1}\,(ST)^{-h_q-2} \, 0 \\
&= T^{-1} \, (ST^{-1})^i \,ST^{-2}\,(ST^{-1})^{h_q} \, ST^{-1}S \, 0 \\
&= T^{-1} \, (ST^{-1})^i \,ST^{-2}\,(ST^{-1})^{h_q} \, 0 \\
&= \lb -1;(-1)^i ,-2, (-1)^{h_q} \rb.
\intertext{Similarly, we have}
(ST)^{h_q+1-i} \, \big(\lb 0;1^{h_q} \rb \big)
&= (ST)^{h_q+1-i} \, (ST)^{h_q} \, 0 \\
&= (ST)^{-i-2} \, 0 \\
&= T^{-1} \, (ST^{-1})^{i} \, ST^{-1}S \, 0 \\
&= T^{-1} \, (ST^{-1})^{i} \, 0 = \lb -1;(-1)^{i} \rb,
\end{align*}
which proves the equations ~\eqref{E4.2a}. The lexicographic order in Section \ref{C3} implies that
\begin{multline*}
\lb -1; \,\rb < \lb -1; -2, (-1)^{h_q} \rb < \lb -1; -1 \rb < \lb -1; -1,-2,(-1)^{h_q} \rb < \ldots \\
\ldots < \lb -1; (-1)^{h_q} \rb < \lb -1; (-1)^{h_q},2,(-1)^{h_q} \rb = - \frac{\lambda_q}{2}.
\end{multline*}
Using the identities \eqref{E4.2a} one can then easily deduce the ordering in \eqref{E4.2b}.
\end{proof}
\begin{lemma}
\label{E5}
For $q=2h_q+3,\,h_q\geq 1$, define the intervals
\begin{align*}
&I_{2i+1}=\big( \lb -1;(-1)^i,-2,(-1)^{h_q}\rb,\frac{\lambda_q}{4}\big) \quad \text{for} \quad 0\leq
i\leq h_q,\\
&I_{2i}=\big(\lb -1;(-1)^i\rb , \frac{\lambda_q}{4} \big)\quad \text{for}\quad 1\leq i\leq h_q \quad \text{and set} \\
&I_{-i}= -I_i \quad \text{for} \quad 1\leq i\leq \kappa_q.
\end{align*}
Then $\Phi_i\subset I_i$ for all $1\leq i\leq \kappa_q$ and
$\vartheta_n(\overline{I}_i)\subset I_j$ for all $n\in \mathcal{N}_{i,j}$ unless $(i,j) = (\pm k,\pm
(k-2))$ with $3\leq k\leq \kappa_q$. In the remaining cases we have that if $3\leq k\leq \kappa_q$ and $n\in \mathcal{N}_{\pm k,\pm (k-2)}$ then
$\vartheta_{n}({I}_{\pm k})\subset I_{\pm (k-2)}$ but $\vartheta_{n}({\overline{I}}_{\pm k})\not\subset I_{\pm (k-2)}$.

\end{lemma}
\begin{proof}
Since the proof of this lemma for the non-exceptional cases proceeds along the same lines as the proof of Lemma \ref{E4.4} (for even $q$),
we only consider the case where $\vartheta_n(\overline{I}_i)\not\subset I_j$ for $n \in \mathcal{N}_{i,j}$.
This happens only for $(i,j)=(\pm k,\pm (k-2))$,
where $\mathcal{N}_{i,j}=\{\pm 1\}$ In all these cases one finds indeed that $\vartheta_1\big(\lb-1;(-1)^i,-2,(-1)^{h_q}\rb\big)=\lb-1;(-1)^{i-1},-2,(-1)^{h_q}\rb$ and
$\vartheta_1\big(\lb-1;(-1)^i\rb\big)=\lb-1;(-1)^{i-1}\rb$. Hence the left boundary point of these intervals is
mapped onto the left boundary point of the image intervals. The case of negative indices $(i,j)$ follows
once more from the symmetry of the intervals and the sets $\mathcal{N}_{i,j}$.
\end{proof}
To finally prove Lemma~\ref{E1.5} one has to enlarge the intervals $I_i$ slightly in order that
$\vartheta_n(\overline{I}_i)\subset I_j$ for all $n\in \mathcal{N}_{i,j}$. In the case $q=2h_q+2$ and
$h_q\geq 2$ one can take the intervals $$I_i=-I_{-i}=\big(\lb -1;(-1)^i,n_i\rb ,\frac{\lambda_q}{4}\big)$$ with
$n_i > n_{i-1}$ for $2\leq i\leq h_q$ and $n_1$ large enough.
 In the case $q=2h_q+3$ and $ h_q\geq 1$ one can choose the intervals
\begin{align*}
&I_{2i+1}=-I_{-2i-1}=\big( \lb -1;(-1)^i,-2,(-1)^{h_q},n_{2i+1}\rb ,\frac{\lambda_q}{4}\big) \quad
\text{for} \quad 0\leq i\leq h_q,\\
&I_{2i}=-I_{-2i}=\big(\lb -1;(-1)^i,n_{2i}\rb , \frac{\lambda_q}{4} \big)\quad \text{for}\quad 1\leq
i\leq h_q.
\end{align*}
with $n_{2i+1} > n_{2i} > n_{2i-1}>n_{2i-2}$ for all $1\leq i\leq h_q$ and $n_1$ large enough.

\medskip
The existence of the discs $D_i$ for all $ i\in A_{\kappa_q}$ in Lemma~\ref{E1.5} shows that the operator
$\mathcal{L}_s$ is well defined on the Banach space $B=\oplus_{i\in A_{\kappa_q}} B(D_i)$ where $B(D_i)$ is the
Banach space of functions which are holomorphic on the disc $D_i$ and continuous on its closure, together with the sup norm.
\begin{theorem}
\label{E6}
The operator $\mathcal{L}_s:B\to B$ is nuclear of order zero for $\re{s}> \frac{1}{2}$ and it extends to a
meromorphic family of nuclear operators of order zero, in the entire complex plane, with poles only at the points
$s_k=\frac{1-k}{2}, \, k=0,1,2 \ldots$.
\end{theorem}
\begin{proof}
It is easy to verify that the operator $\mathcal{L}_s$ can be written as a $2\kappa_q\times 2\kappa_q$ matrix operator which for even $q$ has the form
$$\mathcal{L}_s=
\begin{pmatrix}
0&0&\ldots&0&\mathcal{L}_{2,s}^\infty&\mathcal{L}_{-1,s}^\infty&0&\ldots&0&0\\
\mathcal{L}_{1,s}&0&\ldots&0&\mathcal{L}_{2,s}^\infty&\mathcal{L}_{-1,s}^\infty&0&\ldots&0&0\\
0&\mathcal{L}_{1,s}&0&\ldots&\mathcal{L}_{2,s}^\infty&\mathcal{L}_{-1,s}^\infty&0&\ldots&0&0\\
\vdots&\vdots&\ddots&\vdots& \vdots&\vdots&\vdots&\vdots&\vdots&\vdots\\
0&\ldots&0&\mathcal{L}_{1,s}&\mathcal{L}_{2,s}^\infty&\mathcal{L}_{-1,s}^\infty&0&\ldots&0&0\\
0&0&\ldots&0&\mathcal{L}_{1,s}^\infty&\mathcal{L}_{-2,s}^\infty&\mathcal{L}_{-1,s}&0&\ldots&0\\
\vdots&\vdots&\vdots&\vdots&\vdots&\vdots&\vdots&\ddots&\vdots&\vdots\\
0&0&\ldots&0&\mathcal{L}_{1,s}^\infty&\mathcal{L}_{-2,s}^\infty&\ldots&0&\mathcal{L}_{-1,s}&0\\
0&0&\ldots&0&\mathcal{L}_{1,s}^\infty&\mathcal{L}_{-2,s}^\infty&0&\ldots&0&\mathcal{L}_{-1,s}\\
0&0&\ldots&0&\mathcal{L}_{1,s}^\infty&\mathcal{L}_{-2,s}^\infty&0&\ldots&0&0
\end{pmatrix}.
$$
and for odd $q$
$$\mathcal{L}_s=
\begin{pmatrix}
0&\ldots&0&\mathcal{L}_{2,s}&\mathcal{L}_{3,s}^\infty& \mathcal{L}_{-2,s}^\infty&\mathcal{L}_{-1,s}& 0&\ldots&0\\
0&\ldots&0&0&\mathcal{L}_{2,s}^\infty&\mathcal{L}_{-2,s}^\infty&\mathcal{L}_{-1,s}&0&\ldots&0\\
\mathcal{L}_{1,s}&\ddots&0&0&\mathcal{L}_{2,s}^\infty&\mathcal{L}_{-2,s}^\infty&\mathcal{L}_{-1,s}&0&\ldots&0\\
\vdots&\ddots&\ddots&\vdots&\vdots&\vdots&\vdots&\vdots&\vdots&\vdots\\
0&\ldots&\mathcal{L}_{1,s}&0&\mathcal{L}_{2,s}^\infty&\mathcal{L}_{-2,s}^\infty&\mathcal{L}_{-1,s}&0&\ldots&0\\
0&\ldots&0&\mathcal{L}_{1,s}&\mathcal{L}_{2,s}^\infty&\mathcal{L}_{-2,s}^\infty&0&\mathcal{L}_{-1,s}&\ldots&0\\
\vdots&\vdots&\vdots&\vdots&\vdots&\vdots&\vdots&\ddots&\ddots&\vdots\\
0&\ldots&0&\mathcal{L}_{1,s}&\mathcal{L}_{2,s}^\infty&\mathcal{L}_{-2,s}^\infty&0&0&\ddots&\mathcal{L}_{-1,s}\\
0&\ldots&0&\mathcal{L}_{1,s}&\mathcal{L}_{2,s}^\infty&\mathcal{L}_{-2,s}^\infty&0&0&\ldots&0\\
0&\ldots&0&\mathcal{L}_{1,s}&\mathcal{L}_{2,s}^\infty&\mathcal{L}_{-3,s}^\infty&\mathcal{L}_{-2,s}&0&\ldots&0
\end{pmatrix}
$$
with $\mathcal{L}_{\pm n ,s}^\infty$ and
$\mathcal{L}_{n,s}$ as defined in (\ref{aux}) and (\ref{aux1}).
In a similar manner as for the transfer operator of the Gau{\ss} map (cf.~\cite{Ma90}) one can show that the
operators $\mathcal{L}_{n,s}^\infty ,\ n=\pm 1,\pm 2,\pm 3$ define meromorphic families of nuclear
operators $\mathcal{L}_{n,s}^\infty :B(D_i)\to B(D_j)$ on the Banach spaces of holomorphic functions on the discs $D_i$ and $D_j$ with $\mathcal{N}_{i,j}=\mathbb{Z}_{\geq n}$ (in the case of $+$) and $\mathcal{N}_{i,j}=\mathbb{Z}_{\leq -n}$
(in the case of $-$) for $n=1,2,3$.
These operators have poles only at the points $s = s_k=\frac{1-k}{2}$ for $k=0,1,\ldots$.
Additionally, the operators $\mathcal{L}_{n,s},\, n=\pm 1,\pm 2$ with $\mathcal{L}_{n,s} :B(D_{ i})\to B(D_{j})$ are holomorphic nuclear operators in the entire $s$-plane
on  the corresponding Banach spaces of holomorphic functions on the discs for which $\mathcal{N}_{i,j}=\{\pm n\},\, n=1,2$.
Since a matrix-valued operator with poles in the entries can not have more poles than its respective components, it follows that the operator $\mathcal{L}_s$ has
precisely the desired properties in the Banach space $B=\oplus_{i\in A_{\kappa_q}} B(D_i)$.
\end{proof}

\pagebreak
\section{Reduced transfer operators and functional equations}
\label{F}
\subsection{The symmetry operator $P$}
\label{symm}
From the matrix representation in the proof of Theorem \ref{E6} it can be seen that the transfer operator $\mathcal {L}_s$ possesses a certain symmetry.
In this subsection we will discuss this symmetry in further detail and explain how it can be used to derive functional equations.
For this purpose, define the operator $P:B\to B$ by
\begin{equation*}
\label{P}
 (P\underline{f})_i(z):= f_{-i}(-z) \quad\text{for}\quad \underline{f}= (f_i)_{i\in A_{\kappa_q}}.
\end{equation*}
This operator is well-defined since $D_{-i}=-D_i$ for all $i\in A_{\kappa_q}$. It is also clear that $P^2=id_B$.
That $P$ is indeed a symmetry for the transfer operator follows from the following lemma.
\begin{lemma}
\label{Psym}
The operators $P:B\to B$ and $\mathcal{L}_s:B\to B$ commute for all $s\in \mathbb{C}\setminus \{s_1,s_2,\ldots\}$.
\end{lemma}
\begin{proof}
Let $\re{s} > \frac{1}{2}$ and suppose that $\underline{f}\in B$. To extend $\vartheta_n^\prime$ to the complex discs $D_i$ we use the convention $(n+z)^{2s}:= \big((n+z)^2\big)^s$.
It is then easy to see that
\begin{align*}
&\mathcal{L}_{l,s}(Pf)_i(z)=\sum_{n\ge l} \left(\frac{1}{z+n\lambda_q}\right)^{2s}f_{-i}\left(\frac{1}{z+n\lambda_q}\right)\\
&= \sum_{n\ge l} \left(\frac{1}{-z-n\lambda_q}\right)^{2s}f_{-i}\left(\frac{-1}{-z-n\lambda_q}\right)
   = \mathcal{L}_{-l,s}f_{-i}(-z)
\end{align*}
for any positive integer $l$. One can verify that the matrix elements of $\mathcal{L}_s$ satisfy the identities: $(\mathcal{L}_s)_{i,j}=\mathcal{L}_{l,s}$
if and only if $(\mathcal{L}_s)_{-i,-j}=\mathcal{L}_{-l,s}$ and $(\mathcal{L}_{s})_{i,j}=\mathcal{L}^{\infty}_{l,s}$ if and only if $(\mathcal{L}_s)_{-i,-j}=\mathcal{L}^{\infty}_{-l,s}$.
Combining these two observations, the fact that $\mathcal{L}_s P\underline{f}(z)=P\mathcal{L}_s\underline{f}(z)$ follows immediately.
Since the operators $P\,\mathcal{L}_s$ and $\mathcal{L}_s \,P$ are both meromorphic
in the entire $s$-plane with poles only at the points $s_1,s_2,\ldots$, it follows by meromorphic extension that the identity
 $P\,\mathcal{L}_s = \mathcal{L}_s \,P$ holds for any $s$ in the region $\mathbb{C}\setminus \{s_1,s_2,\ldots\}$.
\end{proof}
The previous lemma allows us to restrict the operator $\mathcal{L}_s$ to the eigenspaces of the operator $P$. Since this is a linear
 involution it can only have eigenvalues $\pm 1$. Denote the corresponding eigenspaces by $B_\pm$. Then
$\underline{f}=(f_i)_{i\in _{A\kappa_q}}\in B_{\pm}$ if and only if $f_{-i}(-z)=\pm \,f_i(z)$ for $ i \in A_{\kappa_q}$.
Let $B_{\kappa_q}$ denote the Banach space $B_{\kappa_q}=\oplus_{1\leq i\leq\kappa_q} B(D_i)$ with
the discs $D_i$ as defined earlier in Lemma~\ref{E1.5}. Then the transfer operator $\mathcal{L}_s$
restricted to the spaces $B_\pm $ induces operators
$\mathcal{L}_{s,\pm} $ on the Banach space $B_{\kappa_q}$.
Let $\overrightarrow{g}=(g_i)_{1\leq i\leq \kappa_q} \in B_{\kappa_q}$.
For $q=2h_q+2$ we get
\begin{align}
\label{reduced1}
&(\mathcal{L}_{s,\pm}\overrightarrow{g})_{1}(z) = L_{2,s}^\infty\, g_{h_q}(z) \pm L_{-1,s}^\infty \,g_{h_q}(z), \nonumber\\
&(\mathcal{L}_{s,\pm}\overrightarrow{g})_{i\hspace{\wonemi}}(z) = L_{1,s}\, g_{i-1}(z)+L_{2,s}^\infty \,
g_{h_q} (z)\pm L_{-1,s}^\infty\, g_{h_q}(z),\, 2\leq i \leq h_q.
\end{align}
For $q=3$ we get
\begin{equation}
\label{reduced2}
(\mathcal{L}_{s,\pm}\overrightarrow{g})_1(z)= L_{3,s}^\infty \, g_1(z)\pm L_{-2,s}^\infty \, g_1(z).
\end{equation}
For $q=2h_q+3>5$ we get
\begin{align}
\label{reduced3}
(\mathcal{L}_{s,\pm}\overrightarrow{g})_1(z)=&L_{2,s}\,g_{2h_q}(z)+ L_{3,s}^\infty
\, g_{\kappa_q}(z) \pm L_{-1,s}\, g_{2h_q}(z)\pm L_{-2,s}^\infty\, g_{\kappa_q}(z), \nonumber\\
(\mathcal{L}_{s,\pm}\overrightarrow{g})_2(z)=&L_{2,s}^\infty\, g_{\kappa_q}(z)\pm
L_{-1,s}\, g_{2h_q}(z)\pm L_{-2,s}^\infty\, g_{\kappa_q}(z), \\
(\mathcal{L}_{s,\pm}\overrightarrow{g})_i(z)=& L_{1,s}\, g_{i-2}(z) +
L_{2,s}^\infty\, g_{\kappa_q}(z)\pm L_{-1,s}\, g_{2h_q}(z)\pm
L_{-2,s}^\infty\, g_{\kappa_q}(z), \nonumber \\
& 3\leq i\leq \kappa_q.\nonumber
\end{align}

For $i>0$ the operators $L_{i,s}^\infty$ and $ L_{i,s}$ coincide with the operators
$\mathcal{L}_{i,s}^\infty$ and $ \mathcal{L}_{i,s}$, whereas $L_{-i,s}^\infty
g(z)=\sum_{n=i}^\infty \frac{1}{(z-n \lambda_q)^{2s}} g\big(\frac{1}{z-n \lambda_q}\big)$ and
$ L_{-i,s}g(z)= \frac{1}{(z-i \lambda_q)^{2s}}g\big(\frac{1}{z-i \lambda_q}\big)$.

\subsection{Functional equations}
It is known that the transfer operator for the Gau{\ss} map induces a
matrix-valued transfer operator for each finite index subgroup, $\Gamma$, of the modular group (cf.~e.g.~\cite{CM00},\cite{CM01}).
This transfer operator can be described by the representation of $\PSL{\IZ}$ induced by the trivial representation of $\Gamma$.
It was shown in \cite{CM01} that eigenfunctions, with eigenvalue $1$, of this induced transfer operator, fulfil vector-valued, finite term functional equations,
analogous to the so-called Lewis equation. This implies that such eigenfunctions are closely related to the period functions of Lewis and Zagier \cite{LZ01} for these groups.

As we will see soon, it is possible to derive similar functional equations for the family of transfer operators considered in this paper.
However, the relationship to period functions, in the case of an arbitrary Hecke triangle group, is not clear at this time.
In the case $q=3$ it was shown in
\cite{BM09} that the solutions of the functional equation derived from our transfer operator
$\mathcal{L}_s$ for $0<\re{s}<1,\, s\neq \frac{1}{2}$ are indeed in a one-to-one correspondence with the Maa{\ss}
cusp forms for $G_3$.

Since the spectrum of the operator $\mathcal{L}_s$ is the union of the
spectra of the two operators $\mathcal{L}_{s,\epsilon},\, \epsilon = \pm 1$, we use these
operators to derive the corresponding functional equations. In the case of $q=3$,
the eigenfunctions, $\overrightarrow{g}=(g_1)$, with eigenvalue $\rho =1$, satisfy the equation
\begin{equation}
\label{eq}
g_1 = g_1|(\mathbb{N}_3+\epsilon \mathbb{N}_{-2})
\end{equation}
where, for $k\ge 1$, we have defined
\begin{align*}
 g_1|\mathbb{N}_{k}(z) & = g_1|\sum_{l=k}^\infty ST^l:= \sum_{l=k}^\infty \left(\frac{1}{z+l
\lambda_q}\right)^{2s}\,
g_1\left(\frac{- 1}{z +l\lambda_q}\right)\, \text{and}\\
g_1|\mathbb{N}_{-k}(z) & = g_1|\sum_{l=k}^\infty \tilde{S}T^{-l}:= \sum_{l=k}^\infty \left(\frac{1}{z-l
\lambda_q}\right)^{2s}\,
g_1\left(\frac{1}{z -l\lambda_q}\right).
\end{align*}
Here $Tz=z+\lambda_q,\, Sz=\frac{-1}{z}$, $Jz=-z$ and $\tilde{S} z = JSz=\frac{1}{z}$.
This action is similar to the usual slash-action of weight $s$ but we have extended it in a natural way to $\mathbb{C}[G_q]$, the extended group ring of $G_q$ over $\mathbb{C}$,
consisting of (possibly countably infinite) formal sums of elements in $G_q$ with coefficients in $\mathbb{C}$.
One now sees that $g_1|\mathbb{N}_3(1-T)=g_1|ST^3$ and $g_1|\mathbb{N}_{-2}(1-T)=-g_1|\tilde{S}T^{-1}$, which leads to the following
four term functional equation
$$ g_1|(1-T)= g_1|(ST^3-\epsilon \tilde{S}T^{-1}).$$
Explicitly, this can be written as
\begin{equation}
\label{eqq}
g_1(z)=g_1(z+1) + \left(\frac{1}{z+3}\right)^{2 s}g_1\left(\frac{-1}{z+3}\right)- \epsilon \left(\frac{1}{z-1}\right)^{2 s}g_1\left(\frac{1}{z-1}\right).
\end{equation}
Using the fact that $JTJ=T^{-1}$ it is easy to verify that every solution of (\ref{eq}) satisfies the equation $g_1(z)=\epsilon g_1(-z-1)$.
Therefore, only solutions $g_1$ of (\ref{eqq}) with this property lead to eigenfunctions of the transfer operator.
On the one hand, it now follows that $g_1$ has to satisfy the shifted four term functional equation which was studied in \cite{BM09}:
\begin{equation}
\label{eqqq}
g_1(z)=g_1(z+1)+\left(\frac{1}{z+3}\right)^{2s}g_1\left(\frac{-1}{z+3}\right)-\left(\frac{1}{z-1}\right)^{2s}g_1\left(\frac{-z}{z-1}\right).
\end{equation}
On the other hand, every solution $g_1$ of (\ref{eqqq}) that additionally satisfies $g_1(z)=\epsilon g_1(-z-1)$ is also a solution of (\ref{eqq}).

For $q=2h_q+2$ and $ h_q\geq 1$ any eigenfunction $\overrightarrow{g}=(g_i)_{1\leq i\leq h_q}$ must satisfy the equation
$g_1=g_{h_q}|(\mathbb{N}_{2} + \epsilon\mathbb{N}_{-1})$. By induction on $i$ it follows that
$$g_i=g_1|P_{i-1}(ST),\, 2\leq i\leq h_q,$$ where $g|P_i(g)$ for $g\in G_q$ is an abbreviation for
$g|P_i(g)=g|\sum_{l=0}^i g^l$.
Hence the function $g_1$ fulfills the equation
$$g_1=g_1|P_{h_q-1}(ST)(\mathbb{N}_{2} + \epsilon\mathbb{N}_{-1}).$$
However, $\mathbb{N}_{2}(1-T)=ST^2$ and $\mathbb{N}_{-1}(1-T)=-\tilde{S}$ in $\mathbb{C}[G_q]$, which leads to the $q$-term
functional equation
\begin{equation*}
 g_1|(1-T)=g_1|P_{h_q-1}(ST)(ST^2-\epsilon \tilde{S}).
\end{equation*}
For $q=4$ this can be written explicitly as
\begin{equation*}
g_1(z)=g_1(z+\lambda_4)+ \left(\frac{1}{z+2 \lambda_4}\right)^{2 s}g_1\left(\frac{-1}{z+2
\lambda_4}\right)- \epsilon \left(\frac{1}{z}\right)^{2 s}g_1\left(\frac{1}{z}\right).
\end{equation*}
For $q=6$ one finds that
\begin{align*}
&g_1(z)=g_1(z+\lambda_6)+ \left(\frac{1}{z+2 \lambda_6}\right)^{2 s}g_1\left(\frac{-1}{z+2
\lambda_6}\right)- \epsilon \left(\frac{1}{z}\right)^{2 s}g_1\left(\frac{1}{z}\right)\\
&+\left(\frac{1}{-\lambda_6 z+1-2 \lambda_6^2}\right)^{2s}g_1\left(\frac{z+2\lambda_6}{-\lambda_6z+1-2 \lambda_6^2}\right)- \epsilon \left(\frac{1}{1+\lambda_6 z}\right)^{2 s}g_1\left(\frac{-z}{1+\lambda_6 z}\right) \nonumber.
\end{align*}
For $q=2h_q+3$ and $ h_q\geq 1$ one finds that
$$
g_1=g_{2h_q}|ST^2+g_{2h_q+1}|\mathbb{N}_{3}+\epsilon g_{2h_q}|\tilde{S}T^{-1}+\epsilon g_{2h_q+1}|\mathbb{N}_{-2}
$$
and
$$
g_2=g_{2h_q+1}|\mathbb{N}_{2}+\epsilon g_{2h_q}|\tilde{S}T^{-1}+\epsilon g_{2h_q+1}|\mathbb{N}_{-2}.
$$
Hence
\begin{equation}
\label{g_1}
g_1=g_2+g_{2h_q}|ST^2-g_{2h_q+1}|ST^2.
\end{equation}
Using induction on $i$, one can show that
$$
g_{2i}=g_2|P_{i-1}(ST),\, 1\leq i\leq h_q
$$
and
$$
g_{2i+1}= g_1|(ST)^i+g_2|P_{i-1}(ST),\, 1\leq i\leq h_q.
$$
In particular, setting $i=h_q$, we can express $g_{2h_g}$ and $g_{2h_q+1}$ in terms of $g_1$ and $g_2$:
\begin{equation*}
g_{2h_q}=g_2|P_{h_q-1}(ST) \quad\text{and}\quad g_{2h_q+1}= g_1|(ST)^{h_q}+g_2|P_{h_q-1}(ST).
\end{equation*}
Inserting these expressions into (\ref{g_1}) leads to an expression for $g_2$ only involving $g_1$:
\begin{equation*}
g_2=g_1+g_1|(ST)^{h_q+1}T.
\end{equation*}
This allows us to express both $g_{2h_q}$ and $g_{2h_q+1}$ in terms of $g_1$:
\begin{equation*}
 g_{2h_q}=g_1 |\left(\id+(ST)^{h_q+1}T \right)P_{h_q-1}(ST)
\end{equation*}
and
\begin{equation*}
 g_{2h_q+1}=g_1|(ST)^{h_q}+\left(g_1+g_1|(ST)^{h_q+1}T \right)P_{h_q-1}(ST).
\end{equation*}
Inserting these expressions into (\ref{g_1}) gives the following functional equation for $g_1$:
\begin{align*}
g_1&=g_1|P_{h_q-1}(ST)ST^2+g_1|(ST)^{h_q+1}T P_{h_q-1}(ST)ST^2+g_1|(ST)^{h_q}\mathbb{N}_{3}\nonumber\\
&+g_1|P_{h_q-1}(ST)\mathbb{N}_{3} +g_1|(ST)^{h_q+1}TP_{h_q-1}(ST) \mathbb{N}_{3}\\
&+\epsilon( g_1|P_{h_q-1}(ST)\tilde{S}T^{-1}+ g_1|(ST)^{h_q+1}T P_{h_q-1}(ST)\tilde{S}T^{-1}+g_1|(ST)^{h_q}
\mathbb{N}_{-2}\nonumber\\
& +g_1|(ST)^{h_q+1}TP_{h_q-1}(ST)\mathbb{N}_{-2}\nonumber+g_1|P_{h_q-1}(ST)\mathbb{N}_{-2}).\nonumber
\end{align*}
Since $\mathbb{N}_{3}(1-T)=ST^3$ and $\mathbb{N}_{-2}=-\tilde{S}T^{-1}$ in $\mathbb{C}[G_q]$, we obtain a "Lewis-type" equation for the Hecke triangle group $G_q$, with $q$ odd, of the following form:
\begin{align*}
g_1|(1-T) &= g_1| (P_{h_q-1}(ST)ST^2+(ST)^{h_q+1}TP_{h_q-1}(ST)ST^2+(ST)^{h_q}ST^3)\nonumber\\
&-\epsilon\, g_1 |(P_{h_q-1}(ST)\tilde{S}+(ST)^{h_q+1}TP_{h_q-1}(ST)\tilde{S}+(ST)^{h_q}\tilde{S}T^{-1}).\nonumber
\end{align*}
For $q=5$ this equation can be written explicitly as
\begin{align*}
g_1|(1-T)&= g_1| ( ST^2+(ST)^2 T ST^2+(ST)^2 T^2)\nonumber\\ &
-\epsilon\, g_1 |(\tilde{S}+(ST)^2 T \tilde{S}+ST\tilde{S}T^{-1}).
\end{align*}

\section{The Selberg zeta function for Hecke triangle groups}
\label{G}
We want to express the Selberg zeta function for the Hecke triangle groups $G_q$ in terms of Fredholm determinants of the transfer operator $\mathcal{L}_s$ for the map $f_q$. Our construction is analogous to that for modular groups and the Gau{\ss} map \cite{CM00}. We start with a discussion of the Ruelle zeta function for the map $f_q$.
\subsection{The Ruelle zeta function and the transfer operator for the map $f_q$}
We have seen that the transfer operator 
for the map $f_q:I_q\to I_q$ can be written
as
\begin{equation}
\label{Lbeta}
(\mathcal{L}_s f)(x) \, = \sum_{n\in \mathbb{Z}:\lb 0;n,x \rb\in\Ar} \big(\vartheta_n^\prime
(x)\big)^s f\big(\vartheta_n(x)\big),
\end{equation}
where $\Ar$ denotes the set of regular $\lambda_q$-CF of all points $x\in
I_q$ and $\lb 0;n,x\rb$ denotes the continued fraction $\lb 0;n,a_1,a_2,\ldots\rb$ for $x=\lb
0;a_1,a_2,\ldots\rb\in I_q$. The iterates $\mathcal{L}_s^k,\, k=1,2,\ldots$ of this operator
have the form
 \begin{equation*}
(\mathcal{L}_s^k f)(x) \, = \sum_{(n_1,\ldots,n_k)\in \mathbb{Z}^k:\lb 0;n_1,\ldots,n_k,x
\rb\in\Ar} \left(\vartheta_{n_1,\ldots,n_k}^\prime
(x)\right)^s f\left(\vartheta_{n_1,\ldots,n_k}(x)\right),
\end{equation*}
where $\vartheta_{n_1,\ldots,n_k}$ denotes the map $\vartheta_{n_1}\circ\ldots\circ \vartheta_{n_k}$.

We have also seen that the set of $k$-tuples $(n_1,\ldots,n_k)\in \mathbb{Z}^k$ with \linebreak
$\lb 0;n_1,\ldots,n_k,x \rb\in\Ar$ only depends on the interval $I_i^{\circ}$ to which  $x$ belongs.
Denote this set by $\mathcal{F}_i^{k}$, i.e.
\begin{equation*}
\mathcal{F}_i^{k}=\{(n_1,\ldots,n_k)\in \mathbb{Z}^k:\lb 0;n_1,\ldots,n_k,x
\rb\in\Ar\quad\text{for all}\quad x \in I_i^\circ \}\,.
\end{equation*}
Hence, for $x\in I_i$, the operator $\mathcal{L}_s^k$ can be written as
\begin{align*}
(\mathcal{L}_s^k f)(x) \, & = \sum_{(n_1,\ldots,n_k)\in\mathcal{F}_i^{k}}
\left(\vartheta_{n_1,\ldots,n_k}^\prime (x) \right)^s f\big(\vartheta_{n_1,\ldots,n_k}(x)\big).
\intertext{Let $f_j$ denote the restriction $f|I_j$ and $\underline{n}_k=(n_1,\ldots,n_k)\in \mathbb{Z}^k$. Then}
(\mathcal{L}_s^k f)_i(x) \, & = \sum_{j\in\mathcal{A}_{\kappa_q}}\sum_{\underline{n}_k\in\mathcal{F}_i^{k}} \left(\vartheta_{\underline{n}_k}^\prime (x) \right)^s \,\chi_{I_j}\left(\vartheta_{\underline{n}_k}(x)\right)\,f_j\left(\vartheta_{\underline{n}_k}(x)\right),
\intertext{which on the Banach space $B=\oplus_{i\in\mathcal{A}_{\kappa_q}} B(D_i)$ reads as}
(\mathcal{L}_s^k \underline{f})_i(z) \, & =
\sum_{j\in\mathcal{A}_{\kappa_q}}\sum_{\underline{n}_k\in\mathcal{F}_i^{k}}
\left(\vartheta_{\underline{n}_k}^\prime (z)\right)^s\, \chi_{D_j}\left(\vartheta_{\underline{n}_k}(z)\right)\,f_j\left(\vartheta_{\underline{n}_k}(z)\right).
\end{align*}
On the Banach space $B$, the trace of $\mathcal{L}_s$ is given by the well-known formula for this type of composition operators (cf.~e.g.~\cite{Ma80})
\begin{equation*}
\text{trace}\, \mathcal{L}_s^k = \sum_{i\in\mathcal{A}_{\kappa_q}}\sum_{\underline{n}_k\in\mathcal{F}_i^{k}}\big(\vartheta_{\underline{n}_k}^\prime (z^\star_{\underline{n}_k}) \big)^s\frac{1}{1-\vartheta_{\underline{n}_k}^\prime (z^\star_{\underline{n}_k})},
\end{equation*}
where $z^\star_{\underline{n}_k}=\lb0;\overline{n_1,\ldots,n_k}\rb$ is the unique (attractive) fixed point of the hyperbolic map
$\vartheta_{n_1,\ldots,n_k}:D_i\to D_i$. However, since these fixed points are in one-to-one correspondence
with the periodic points of period $k$ of the map $f_q$, the following identity holds
\begin{align}
\label{Ruelleze}
\text{trace}\, \mathcal{L}_s^k -\text{trace}\, \mathcal{L}_{s+1}^k & =
\sum_{i\in\mathcal{A}_{\kappa_q}}\sum_{\underline{n}_k\in\mathcal{F}_i^{k}}\left((\vartheta_{n_1}\circ\ldots\
\circ \vartheta_{n_k})^\prime (z^\star_{\underline{n}_k}) \right)^s
\nonumber\\
&=\sum_{i\in\mathcal{A}_{\kappa_q}}\sum_{\underline{n}_k\in\mathcal{F}_i^{k}}\prod_{l=1}^k\left(\vartheta_{n_l}^\prime(\vartheta_{n_{l+1}}\circ\ldots\
\circ \vartheta_{n_k}(z^\star_{\underline{n}_k})\right)^s\nonumber\\
&=\sum_{z^\star\in \text{Fix} f_q^k}\prod_{l=0}^{k-1}\big( f_q^\prime\left(f_q^l(z^\star)\right)\big)^{-s}.
\end{align}
We are now able to relate the Ruelle zeta function to the transfer operator.
\begin{proposition}
The Ruelle zeta function, $\zeta_R(s)=\exp \sum_{n=1}^\infty \frac{1}{n} Z_n(s)$, for the Hurwitz-Nakada map $f_q$ can be written in the entire complex $s$-plane as
$$\zeta_R(s)=\frac{\det (1-\mathcal{L}_{s +1})}{\det (1-\mathcal{L}_{s })}$$ with $\mathcal{L}_s :B\to B$ defined in Theorem \ref{E6}.
\end{proposition}
\begin{proof}
Comparing equations (\ref{Ruelleze}) and (\ref{Ruellez}) we see that for $\re{s} $ large enough we have $Z_n(s)= \text{trace} \, \mathcal{L}_s^n -\text{trace} \,
\mathcal{L}_{s+1}^n $ and therefore $\zeta_R(s)=\frac{\det (1-\mathcal{L}_{s +1})}{\det (1-\mathcal{L}_{s })}$. Since the operators $\mathcal{L}_{s}$ are meromorphic and nuclear in the entire $s$-plane,
both Fredholm determinants also have meromorphic continuations. This proves the proposition.
\end{proof}
\subsection{The transfer operator $\mathcal{K}_s$}
\label{secK}
As we have seen, there is a one-to-one correspondence between the closed orbits of the map
$f_q$ and the closed orbits of the geodesic flow, except for the contributions corresponding to the two points $r_q$ and $-r_q$.
These points are not equivalent under the map $f_q$ but they are under the group $G_q$, which means that they correspond to the same closed orbit of the
geodesic flow. We also saw that the Ruelle zeta function for the map $f_q$ contains contributions of both orbits while the Selberg zeta function
 \eqref{SZ-one-orbit} does not.
Thus, to relate the Ruelle and the Selberg zeta functions we have to subtract the contribution of one of these two orbits.
To be specific, we subtract the contribution of the orbit $\mathcal{O}_+$, given by the point $r_q$,
from all the partition functions $Z_{l\kappa_q}(s)$, $l=1,2,\ldots$.
Consider therefore the corresponding Ruelle zeta function, $\zeta^\mathcal{O_+}_R(s)= \exp \left( \sum_{n=1}^\infty \frac{1}{n} Z_n^{\mathcal{O}_+}(s)\right)$, with
\begin{equation}
\label{partO}
Z_n^\mathcal{O}(s)=
\begin{cases}
0 &\text{if $\kappa_q\nmid$ n,}\\
\sum_{x\in {\mathcal{O}_+}}\exp \left(-s \sum_{k=0}^{n-1}\ln f_q^\prime(f_q^k(x))\right)& \text{if $\kappa_q\mid n$}.
\end{cases}
\end{equation}
If $n=l \kappa_q$ we find that $Z_{l \kappa_q}^{\mathcal{O}_+}(s)= \kappa_q \exp (-s l r_{{\mathcal{O}_+}})$
and hence $\zeta^{\mathcal{O}_+}_R(s)=\frac{1}{1-\exp (-s r_{{\mathcal{O}_+}})}$ where $r_{{\mathcal{O}_+}}=
\ln (f_q^{\kappa_q})^\prime (r_q)$.
We now define a transfer operator, $\mathcal{L}_s^{\mathcal{O}_+}:B_{\kappa_q}\to B_{\kappa_q}$, with $B_{\kappa_q}=\oplus_{i=1}^{\kappa_q} B(D_i)$ as in Section \ref{symm},
corresponding to $\mathcal{O}_{+}$.
For $q=2h_q+2$, on the one hand, we set
\begin{align}
\label{transfereven}
(\mathcal{L}_s^{{\mathcal{O}_+}}\overrightarrow{g})_i(z)&=\mathcal{L}_{1,s}g_{i+1}(z),\, 1\leq i\leq
h_q-1,\nonumber\\
(\mathcal{L}_s^{\mathcal{O}_+}\overrightarrow{g})_{h_q}(z)&=\mathcal{L}_{2,s}g_{1}(z).
\end{align}
For $q=2h_q+3$, on the other hand, we set
\begin{align}
\label{transferodd}
(\mathcal{L}_s^{{\mathcal{O}_+}}\overrightarrow{g})_i(z)=& \mathcal{L}_{1,s}g_{i+1}(z),\, 1\leq i \leq h_q,\nonumber\\
(\mathcal{L}_s^{{\mathcal{O}_+}}\overrightarrow{g})_{h_q+1}(z)=&\mathcal{L}_{2,s}g_{h_q+2}(z),\nonumber\\
(\mathcal{L}_s^{{\mathcal{O}_+}}\overrightarrow{g})_{h_q+i}(z)=&\mathcal{L}_{1,s}g_{h_q+i+1}(z),\, 2\leq i\leq h_q,\\
(\mathcal{L}_s^{{\mathcal{O}_+}}\overrightarrow{g})_{2h_q+1}(z)=&\mathcal{L}_{2,s}g_{1}(z).\nonumber
\end{align}
In the first case, the operator $\mathcal{L}_s^{{\mathcal{O}_+}}$ has the form

$$\mathcal{L}_s^{\mathcal{O}_+}=
\begin{pmatrix}
0&\mathcal{L}_{1,s}&0&\ldots&0&0\\
0&0&\mathcal{L}_{1,s}&\ddots&0&0\\
\vdots&\ddots&\ddots&\ddots&\ddots&\vdots&\\
\vdots&\ddots&\ddots&\ddots&\mathcal{L}_{1,s}&0\\
0& 0&0&\ddots&0&\mathcal{L}_{1,s}\\
\mathcal{L}_{2,s}&0&0&\ldots&0&0\\
\end{pmatrix},
$$
while in the second case
$$\mathcal{L}_s^{\mathcal{O}_+}=
\begin{pmatrix}
0&\mathcal{L}_{1,s}&0&\ldots&0&\ldots&0&0&0\\
\vdots&\ddots&\ddots&\ddots&\ddots&\vdots&\vdots&\vdots&\vdots\\
0&0&\ddots&\mathcal{L}_{1,s}&0&\ldots&0&0&0\\
0&0&\ldots&0&\mathcal{L}_{2,s}&0&\ldots&0&0\\
0& 0&0&\ddots&0&\mathcal{L}_{1,s}&0&\ldots&0\\
\vdots&\vdots&\vdots&\vdots&\vdots&\ddots&\ddots&\ddots&\vdots\\
\vdots&\vdots&\vdots&\vdots&\vdots&\ddots&\ddots&\mathcal{L}_{1,s}&0\\
0&\ldots&\ldots&\ldots&\ldots&\ldots&\ddots&0&\mathcal{L}_{1,s}\\
\mathcal{L}_{2,s}&0&\ldots&\ldots&\ldots&\ldots&\ldots&0&0\\
\end{pmatrix}.
$$
To show that these operators indeed correspond to the contribution of the orbit $\mathcal{O}_{+}$ we need to compute the traces of their iterates.
This is achieved by the following
\begin{lemma}
\label{traceO}
The trace of the operator $(\mathcal{L}_s^{\mathcal{O}_+})^n$ is given by
\begin{equation*}
\text{trace}\,(\mathcal{L}_s^{\mathcal{O}_+})^n =
\begin{cases} 0 \quad &\text{for}\quad \kappa_q \nmid n,\\
\kappa_q \,\text{trace}\,\big(\mathcal{L}_{1,s}^{h_q-1}\mathcal{L}_{2,s}\big)^l\quad &\text{for}\quad n=l\kappa_q,
\end{cases}
\end{equation*}
for $q=2h_q+2$ ($\kappa_q=h_q$) and by
\begin{equation*}
\text{trace}\,(\mathcal{L}_s^{\mathcal{O}_+})^n =
\begin{cases} 0 \quad &\text{for}\quad \kappa_q \nmid n,\\
\kappa_q \,\text{trace}\,\big(\mathcal{L}_{1,s}^{h_q}\mathcal{L}_{2,s}\mathcal{L}_{1,s}^{h_q-1}\mathcal{L}_{2,s}\big)^l\quad &\text{for}\quad n=l\kappa_q,
\end{cases}
\end{equation*}
for $q=2h_q+3$ ($\kappa_q=2h_q+1$).
\end{lemma}
\begin{proof}
Since the proof for odd $q$ is completely analogous to that for even $q$, we restrict ourselves to the case of $q=2h_q+2$. Using induction on $i$ it can be shown that
for $1\leq j \leq h_q$ and $\overrightarrow{g}=(g)_{1\leq j\leq h_q}$ one has
\begin{align*}
((\mathcal{L}_s^{\mathcal{O}_+})^i\overrightarrow{g})_j& = \mathcal{L}_{1,s}^i\, g_{i+j},\quad 1\leq j\leq h_q-i,\nonumber\\
((\mathcal{L}_s^{\mathcal{O}_+})^i\overrightarrow{g})_j& = \mathcal{L}_{1,s}^{h_q-j}\,\mathcal{L}_{2,s}\,\mathcal{L}_{1,s}^{i+j-h_q-1}\, g_{i+j-h_q}, \quad
h_q-i+1\leq j \leq h_q.
\end{align*}
This also shows that
\begin{align*}
((\mathcal{L}_s^{\mathcal{O}_+})^{h_q}\overrightarrow{g})_1& = \mathcal{L}_{1,s}^{h_q-1}\mathcal{L}_{2,s}\,g_1,
\nonumber \\
((\mathcal{L}_s^{\mathcal{O}_+})^{h_q}\overrightarrow{g})_j& = \mathcal{L}_{1,s}^{h_q-j}\,\mathcal{L}_{2,s}\,\mathcal{L}_{1,s}^{j-1}\,g_j,\quad 2\leq j\leq h_q,
\end{align*}
and therefore
\begin{align*}
\text{trace}\,(\mathcal{L}_s^{\mathcal{O}_+})^n &= 0 \quad\text{if}\quad h_q\nmid n, \nonumber\\
\text{trace}\,(\mathcal{L}_s^{\mathcal{O}_+})^n &= h_q \,\text{trace}\, (\mathcal{L}_{1,s}^{h_q-1}\,\mathcal{L}_{2,s})^l \quad\text{if} \quad n=l\, h_q= l\, \kappa_q.
\end{align*}
This proves the Lemma.
\end{proof}
Since $\text{trace}\, (\mathcal{L}_{1,s}^{h_q-1}\,\mathcal{L}_{2,s})^l= \text{trace}\,(\mathcal{L}_{2,s}\, \mathcal{L}_{1,s}^{h_q-1})^l$ and
\begin{equation*}
(\big(\mathcal{L}_{2,s}\, \mathcal{L}_{1,s}^{h_q-1})^l\,g\big)(z) = \left(\frac{d}{dz}(\vartheta_{1^{h-1},2})^l\,(z)\right)^s \, g\left((\vartheta_{1^{h-1},2})^l(z)\right)
\end{equation*}
(recall that $\vartheta_{1^{h-1},2}=\vartheta_1\circá¸·\ldots\circ \vartheta_1\circ\vartheta_2$) we see that
\begin{equation*}
\text{trace}\,(\mathcal{L}_{2,s}\, \mathcal{L}_{1,s}^{h_q-1})^l=\left(\frac{d}{dz}(\vartheta_{1^{h-1},2})^l\,(z^\star)\right)^s \frac{1}{1-\frac{d}{dz}(\vartheta_{1^{h-1},2})^l(z^\star)},
\end{equation*}
where $z^\star$ is the attractive fixed point of $\vartheta_1\circá¸·\ldots\circ \vartheta_1\circ\vartheta_2$, i.e. $z^\star = r_q$. Thus
\begin{align*}
\text{trace}\,\big(\mathcal{L}_{2,s}\, \mathcal{L}_{1,s}^{h_q-1}\big)^l-\text{trace}\,\big(\mathcal{L}_{2,s+1}\, \mathcal{L}_{1,s+1}^{h_q-1}\big)^l
 &=\left(\frac{d}{dz}(\vartheta_{1^{h-1},2})^l\,(z^\star)\right)^s \\
 &= \left(\frac{d}{dz}(\vartheta_{1^{h-1},2})\,(z^\star)\right)^{l\,s}.
\end{align*}

\begin{lemma}
The partition function $Z_n^{{\mathcal{O}_+}}$ in (\ref{partO}) can be expressed in terms of the transfer operators
$\mathcal{L}_s^{{\mathcal{O}_+}}$ in (\ref{transfereven}) and (\ref{transferodd}) as
$$Z_n^{{\mathcal{O}_+}}(s)=\text{trace}\,(\mathcal{L}_s^{\mathcal{O}_+})^n -\text{trace}\, (\mathcal{L}_{s +1}^{\mathcal{O}_+})^n. $$
\end{lemma}
\begin{proof}
Without loss of generality, we once again restrict ourselves to the case of $q=2h_q+2$. Since $(f_q^{\kappa_q})^\prime (z^\star)= \big((\vartheta_1\circá¸·\ldots\circ \vartheta_1\circ\vartheta_2)^\prime (z^\star)\big)^{-1}$
and $Z_n^{{\mathcal{O}_+}}= \kappa_q \exp(-s l \, r_{{\mathcal{O}_+}})$, where
$r_{{\mathcal{O}_+}}=\sum_{k=0}^{\kappa_q-1}\ln f_q^\prime(f_q^k (z^\star))= \ln \prod_{k=0}^{\kappa_q-1}f_q^\prime(f_q^k ( z^\star ))=
\ln (f_q^{\kappa_q})^\prime(z^\star)$, we find that
\begin{align*}
Z_n^{{\mathcal{O}_+}} &= \kappa_q \left((f_q^{\kappa_q})^\prime (z^\star)\right)^{-l\,s}=\kappa_q \big((\vartheta_1\circá¸·\ldots\circ \vartheta_1\circ\vartheta_2)^\prime(z^\star)\big)^{l\,s}\\
&= \kappa_q \,\left(\text{trace}\,\big(\mathcal{L}_{2,s}\, \mathcal{L}_{1,s}^{h_q-1}\big)^l-\text{trace}\,\big(\mathcal{L}_{2,s+1}\, \mathcal{L}_{1,s+1}^{h_q-1}\big)^l\right)\nonumber\\
&= \text{trace}\,\big(\mathcal{L}_s^{\mathcal{O}_+}\big)^n -\text{trace}\, \big(\mathcal{L}_{s +1}^{\mathcal{O}_+}\big)^n.
\end{align*}
\end{proof}
It follows that the Ruelle zeta function for ${\mathcal{O}_+}$ can be expressed as
$$
\zeta^{\mathcal{O}_+}_R(s)= \frac{\text{det}\, (1-\mathcal{L}_{s+1}^{\mathcal{O}_+})}{\text{det}\, (1-\mathcal{L}_{s}^{\mathcal{O}_+})}.
$$
We can furthermore show
\begin{lemma}
\label{lemma:det_op}
The Fredholm determinant $\text{det}\, (1-\mathcal{L}_{s}^{\mathcal{O}_+})$ coincides with the Fredholm
determinant $\text{det}\, (1-\mathcal{L}_{1,s}^{h_q-1}\,\mathcal{L}_{2,s} )$ for even $q$ and with
 $\text{det}\, (1-\mathcal{L}_{1,s}^{h_q}\,\mathcal{L}_{2,s}\,\mathcal{L}_{1,s}^{h_q-1} \,\mathcal{L}_{2,s} )$ for odd $q\ge 5$.
Note that $r_q$ is given in these two cases by $\lb0;\overline{1^{h_q-1},2}\rb$ and $\lb0;\overline{1^{h_q},2,1^{h_q-1},2}\rb$, respectively.
\end{lemma}
\begin{proof}
This lemma follows immediately from Lemma \ref{traceO} and the formula \linebreak
$-\ln \det(1-L)=\sum_{n=1}^\infty \frac{1}{n} \text{trace} L^n$, which holds for any trace class operator $L$.
\end{proof}
\begin{remark}
The spectra of the two operators $\mathcal{L}_{s}^{\mathcal{O}_+}$ and $\mathcal{L}_{1,s}^{h_q-1}\,\mathcal{L}_{2,s}$ can be related in a simple fashion as follows: definition (\ref{transfereven}) of the operator $\mathcal{L}_{s}^{\mathcal{O}_+}$ implies that
any eigenfunction $\overrightarrow{g}=(g_i)_{1\leq i\leq h_q}$ with eigenvalue $\rho$ of $\mathcal{L}_{s}^{\mathcal{O}_+}$ fulfills the equation
$\rho^{h_q-1}g_1= \mathcal{L}_{1,s}^{h_q-1}g_{h_q}$ and hence also $\rho^{h_q}g_1= \mathcal{L}_{1,s}^{h_q-1}\,\mathcal{L}_{2,s} g_1$.
Therefore, on the one hand, any eigenvalue $\rho$ of the operator $\mathcal{L}_{s}^{\mathcal{O}_+}$ determines an eigenvalue $\rho^{h_q}$ of the
operator $\mathcal{L}_{1,s}^{h_q-1}\,\mathcal{L}_{2,s}$.
On the other hand, given an eigenfunction $g$ of $\mathcal{L}_{1,s}^{h_q-1}\,\mathcal{L}_{2,s}$ with eigenvalue $\rho= |\rho|\exp(i\alpha)$,
let $\overrightarrow{g}^{(j)}\in B_{\kappa_q}$ be defined by $g^{(j)}_1=g$ and $g^{(j)}_i=\rho_j^{-(h_q+1-i)}\mathcal{L}_{1,s}^{h_q-i}\,\mathcal{L}_{2,s}\,g$ for $ 2\leq i\leq h_q$.
Then $\overrightarrow{g}^{(j)}$ is an eigenfunction of the operator $\mathcal{L}_{s}^{\mathcal{O}_+}$ with eigenvalue $\rho_j$, where $\rho_1,\ldots,\rho_{h_q}$ are the
$h_q$-th roots of $\rho$. This shows that the numbers
$$\rho_j=\sqrt[h_q]{\mid\rho\mid}\exp\left(i \frac{\alpha}{h_q}\right) \
\exp \left(2\pi i\frac{j}{h_q}\right),\, 0\leq j \leq h_q-1$$
are eigenvalues of this operator. This argument provides an alternative proof of the fact that $\text{det}\, (1-\mathcal{L}_{s}^{\mathcal{O}_+})=\text{det}\, (1-\mathcal{L}_{1,s}^{h_q-1}\,\mathcal{L}_{2,s} )$, i.e., a proof of Lemma \ref{lemma:det_op} in the case of even $q$. The corresponding ``direct proof'' for the case of odd $q$ is analogous.
\end{remark}
From the previous lemma it is clear that the contribution of the periodic orbit of the geodesic flow corresponding to ${\mathcal{O}_+}$,
which appears twice in the Fredholm determinant of $\mathcal{L}_s$, is given
by $\text{det}\,(1-\mathcal{L}_{1,s}^{h_q-1}\,\mathcal{L}_{2,s} )$ and $\text{det}(1-\mathcal{L}_{1,s}^{h_q}\,\mathcal{L}_{2,s}\,\mathcal{L}_{1,s}^{h_q-1} \,\mathcal{L}_{2,s} )$ for even
and odd $q\ge 5$, respectively. We thus arrive at the following theorem.
\begin{theorem}
\label{main-theorem}
The Selberg zeta function $Z_S(s)$ for the Hecke triangle group $G_q$ can be written as
$$Z_S(s)=\frac{\text{det}\,(1-\mathcal{L}_s)}{\text{det}\,(1-\mathcal{K}_s)}=\frac{\text{det}\,[(1-\mathcal{L}_{s,+})(1-\mathcal{L}_{s,-})]}{\text{det}\,(1-\mathcal{K}_s)},$$
where $\mathcal{L}_s$, $\mathcal{L}_{s,\pm}$ and $\mathcal{K}_s=\mathcal{L}_{s}^{\mathcal{O}_+}$ are the transfer operators given by Theorem \ref{E6}, (\ref{reduced1})--(\ref{reduced3}) and
(\ref{transfereven})--(\ref{transferodd}), respectively.
\end{theorem}
\begin{proposition}
The spectrum of $\mathcal{K}_s$ is given by $$\sigma(\mathcal{K}_s)=\left\{\prod_{l=0}^{\kappa_q-1}\big(f_q^l(r_q)\big)^{2s+2n},\, n=0,1,2,\ldots\right\} $$
where $\kappa_q$ denotes the period of the point $r_q$
\end{proposition}

\begin{proof}
The spectrum $\sigma(L)$ of a composition operator of the general form $Lf(z)=\varphi(z) f\big(\psi(z)\big)$ on
a Banach space $B(D)$ of holomorphic functions on a domain $D$ with $\psi(\overline{D})\subset D$ is
given by $\sigma(L)=\{\varphi(z^\star) \psi^\prime(z^\star)^n,\,
n=0,1,\ldots\}$ (cf. e.g.~\cite{Ma80}) where $z^\star$ is the unique fixed point of $\psi$ in $D$. For $q=2h_q+2$ the operator
$\mathcal{K}_s$ has this form with $\psi(z)=\vartheta_2\circ\vartheta_1^{h_q-1}(z)$ and $\varphi(z)=\big(\psi^\prime(z)\big)^s$.
Therefore $z^\star=\lb 0;\overline{2,1^{h_q-1}}\rb$ and $(\vartheta_2\circ\vartheta_1^{h_q-1})^\prime(z^\star)=
\vartheta_2^\prime\big(\vartheta_1^{h_q-1}(z^\star)\big)\prod_{l=1}^{h_q-1}\vartheta_1^\prime \big((\vartheta_1)^{h_q-1-l}(z^\star)\big)$.
Since $\vartheta_m^\prime \big(\vartheta_{m}^{-1}(z)\big)=z^2$ for any $ m\in\mathbb{N}$,
$\vartheta_1^{h_q-1}(z^\star)=\vartheta_{2}^{-1}(z^\star)$ and $(\vartheta_1)^{h_q-1-l}(z^\star)=(\vartheta_1)^{-1}\vartheta_1^{h_q-l}(z^\star)$ it follows immediately that
$$(\vartheta_2\circ\vartheta_1^{h_q-1})^\prime(z^\star)=(z^\star)^2\prod_{l=1}^{h_q-1}\big(\vartheta_1^{h_q-l}(z^\star)\big)^2=\prod_{l=0}^{h_q-1}\big(f_q^l(z^\star)\big)^2=\prod_{l=0}^{h_q-1}\big(f_q^l(r_q)\big)^2.$$

\end{proof}
\begin{remark}
Using the explicit form of the maps which fix $r_q$, cf.~e.g.~\cite[Rem.~27]{MS08} (where the upper right entry of the matrix for even $q$ should read $\lambda-\lambda^3$) one can prove that the
spectrum of the operator $\mathcal{K}_s$ can also be written as \linebreak
 \[\left\{\mu_{n}=l^{2s+2n},\,n=0,1,\ldots\right\},\] where
\begin{align*}
l &=\frac{\sqrt{4-\lambda_q^2}}{R\lambda_q+2}=\sqrt{\frac{2-\lambda_q}{2+\lambda_q}}\quad \text{for even $q$ ,}\\
l &=\frac{2-\lambda_q}{R\lambda_q+2}=\frac{2-\lambda_q}{2+R\lambda_q}\quad \text{for odd $q$}.
\end{align*}
\end{remark}
The Selberg zeta function $Z_S(s)$ for Hecke triangle groups $G_q$ and small $q$ has been calculated numerically using the transfer operator $\mathcal{L}_s$ in \cite{St08}.
Besides the case $q=3$, that is, the modular group $G_3=\PSL{\IZ}$ \cite{BM09}, we do not yet know how the
eigenfunctions of the
transfer operator $\mathcal{L}_s$ with eigenvalue $\rho = 1$ are related to the automorphic functions for a general Hecke group $G_q$.
But since the divisor of $Z_S(s)$ is closely related to the automorphic forms on $G_q$ (see for instance \cite[p.~498]{H83})
one would expect that there exist explicit relationships
also for $q>3$, similar to those obtained for the modular group, i.e.  between eigenfunctions of the transfer
operator $\mathcal{L}_{s}$ with eigenvalue one and automorphic forms related to the divisors of $Z_S$ at these $s$-values.

Another interesting problem one could study is the behaviour of the transfer operator
$\mathcal{L}_s$ in the limit when $q$ tends to $\infty$. In this limit the Hecke triangle group $G_q$ tends to
the theta group $\Gamma_{\theta}$, generated by $Sz=\frac{-1}{z}$ and $Tz=z+2$. This group is conjugate
to the Hecke congruence subgroup $\Gamma_0(2)$, for which a
 transfer operator has been constructed in \cite{HMM05} and \cite{FMM07}. One should understand how these two different transfer operators are related to each other.
The limit $q\to\infty$ is quite singular,
since the group $\Gamma_\theta$ has two cusps whereas all the Hecke triangle groups have only one
cusp. Therefore one expects in this limit the appearance of the singular behavior which Selberg predicted already in
\cite{Se89}. Understanding the limit $q\to \infty$ could also shed new light on the Phillips-Sarnak
conjecture \cite{PS85} on the existence of Maa{\ss} wave forms for general non-arithmetic Fuchsian groups.

\bibliographystyle{amsalpha}

\begin{thebibliography}{99999}

\raggedright

\bibitem{BLZ09}
R.W.\ Bruggeman, J.\ Lewis and D.\ Zagier,
\newblock \textit{Period functions for Maa{\ss} wave forms. II: Cohomology},
\newblock in preparation.

\bibitem{BM09}
R.W.\ Bruggeman and R.\ M\"uhlenbruch,
\newblock \textit{Eigenfunctions of transfer operators and cohomology}.
\newblock Journal of Number Theory \textbf{129} (2009), 158--181. \newline
\newblock \href{http://www.ams.org/mathscinet-getitem?mr=2468476}{\texttt{\small MR2468476}},
          \href{http://dx.doi.org/10.1016/j.jnt.2008.08.003}{\texttt{\small doi:10.1016/j.jnt.2008.08.003}}

\bibitem{CM00}
C.H.\ Chang and D.\ Mayer,
\newblock \textit{Thermodynamic formalism and Selberg's zeta function for modular groups}.
\newblock Regul.\ Chaotic Dyn. \textbf{5} (2000), no.3, 281--312. \newline
\newblock \href{http://www.ams.org/mathscinet-getitem?mr=1789478}{\texttt{\small MR1789478}},
          \href{http://dx.doi.org/10.1070/rd2000v005n03ABEH000150}{\texttt{\small doi:10.1070/rd2000v005n03ABEH000150}}

\bibitem{CM01}
C.H.\ Chang and D.\ Mayer,
\newblock \textit{Eigenfunctions of the transfer operators and the period functions for modular groups}.
\newblock in \textit{Dynamical, spectral, and arithmetic zeta functions} (San Antonio, TX, 1999).
\newblock Contemp.\ Math.\ \textbf{290} (2001), 1--40. \newline
\newblock \href{http://www.ams.org/mathscinet-getitem?mr=1868466}{\texttt{\small MR1868466}}

\bibitem{FMM07}
M.\ Fraczek, D.\ Mayer and T.\ M\"uhlenbruch,
\newblock \textit{A realization of the Hecke algebra on the space of period functions for $\Gamma_0(n)$}.
\newblock J.\ reine angew.\ Math.\ \textbf{603} (2007), 133--163. \newline
\newblock \href{http://www.ams.org/mathscinet-getitem?mr=2312556}{\texttt{\small MR2312556}},
          \href{http://dx.doi.org/10.1515/CRELLE.2007.014}{\texttt{\small doi:10.1515/CRELLE.2007.014}}

\bibitem{H83}
D.\ Hejhal,
\newblock \textit{The Selberg trace formula for $\PSL{\IZ}$, Vol. 2}.
\newblock Springer Lecture Notes in Mathematics \textbf{1001},
\newblock Springer Verlag, 1983.\newline
\newblock \href{http://www.ams.org/mathscinet-getitem?mr=0711197}{\texttt{\small MR0711197}}

\bibitem{HMM05}
J.\ Hilgert, D. \ Mayer and H.\ Movasati,
\newblock \textit{Transfer operators for $\Gamma_0(n)$ and the Hecke operators for the period functions of $\PSL{\IZ}$}.
\newblock Math.\ Proc.\ Camb.\ Phil.\ Soc.\ \textbf{139} (2005), 81--116. \newline
\newblock \href{http://www.ams.org/mathscinet-getitem?mr=2155506}{\texttt{\small MR2155506}},
          \href{http://dx.doi.org/10.1017/S0305004105008480}{\texttt{\small doi:10.1017/S0305004105008480}}

\bibitem{Hu89}
A.\ Hurwitz,
\newblock \textit{U\"ber eine besondere Art der Kettenbruch-Entwickelung reeller Gro\"ssen}
\newblock Acta Math.\ \textbf{12} (1889), 367--405. \newline
\newblock \href{http://dx.doi.org/10.1007/BF02391885}{\texttt{\small doi:10.1007/BF02391885}}

\bibitem{LZ01}
J.\ Lewis and D.\ Zagier,
\newblock \textit{Period functions for Maa{\ss} wave forms.I.}
\newblock Ann.\ of Math.\ \textbf{153} (2001), 191--258. \newline
\newblock \href{http://www.ams.org/mathscinet-getitem?mr=1826413}{\texttt{\small MR1826413}},
          \href{http://dx.doi.org/10.2307/2661374}{\texttt{\small doi:10.2307/2661374}}

\bibitem{Ma90}
D.\ Mayer,
\newblock \textit{On the thermodynamic formalism for the {G}auss map}.
\newblock Comm.\ Math.\ Phys, \textbf{130}, No 2, (1990) , 311-333. \newline
\newblock \href{http://www.ams.org/mathscinet-getitem?mr=1059321}{\texttt{\small MR1059321}},
          \href{http://projecteuclid.org/euclid.cmp/1104200514}{\texttt{\small euclid.cmp/1104200514}}

\bibitem{Ma80}
D.\ Mayer,
\newblock \textit{On composition operators on Banach spaces of holomorphic functions}.
\newblock Journal of Functional Analysis \textbf{35}, No 2, (1980) , 191--206. \newline
\newblock \href{http://www.ams.org/mathscinet-getitem?mr=0561985}{\texttt{\small MR0561985}},
          \href{http://dx.doi.org/10.1016/0022-1236(80)90004-X}{\texttt{\small doi:10.1016/0022-1236(80)90004-X}}

\bibitem{MM09}
D.\ Mayer and T.\ M\"uhlenbruch,
\newblock \textit{Nearest $\lambda_q$-multiple fractions}.
\newblock To appear in: Proceedings of the Spectrum and Dynamics workshop, CRM Proceedings and Lecture Notes Series, (CRM and AMS), Montreal 2008. \newline
\newblock \href{http://arxiv.org/abs/0902.3953}{\texttt{\small arXiv:0902.3953}}

\bibitem{MS08}
D.\ Mayer and F.\ Str\"omberg,
\newblock \textit{Symbolic dynamics for the geodesic flow on Hecke surfaces}.
\newblock Journal of Modern Dynamics \textbf{2} (2008), 581--627. \newline
\newblock \href{http://www.ams.org/mathscinet-getitem?mr=2449139}{\texttt{\small MR2449139}},
          \href{http://dx.doi.org/10.3934/jmd.2008.2.581}{\texttt{\small doi:10.3934/jmd.2008.2.581}}

\bibitem{Na95}
H.\ Nakada,
\newblock \textit{Continued fractions, geodesic flows and Ford circles},
\newblock in \textit{Algorithms, Fractals, and Dynamics}, 179--191, Edited by T.\ Takahashi, Plenum Press, New York, 1995. \newline
\newblock \href{http://www.ams.org/mathscinet-getitem?mr=1402490}{\texttt{\small MR1402490}}

\bibitem{PS85}
R.\ Phillips and P.\ Sarnak,
\newblock \textit{On cusp forms for cofinite subgroups of $PSL(2,\mathbb{R})$},
\newblock Invent.\ Math.\ \textbf{80} (1985), 339--364.\newline
\newblock \href{http://www.ams.org/mathscinet-getitem?mr=0788414}{\texttt{\small MR0788414}},
          \href{http://dx.doi.org/10.1007/BF01388610}{\texttt{\small doi:10.1007/BF01388610}}

\bibitem{Ro54}
D.\ Rosen,
\newblock \textit{A class of continued fractions associated with certain properly discontinuous groups},
\newblock Duke Math.\ J.\ \textbf{21} (1954) , 549--563. \newline
\newblock \href{http://www.ams.org/mathscinet-getitem?mr=MR0065632}{\texttt{\small MR0065632}},
          \href{http://projecteuclid.org/euclid.dmj/1077465884}{\texttt{\small euclid.dmj/1077465884}}

\bibitem{RS92}
D.\ Rosen and T. A.\ Schmidt,
\newblock \textit{Hecke groups and continued fractions},
\newblock Bull.\ Austral.\ Math.\ Soc.\ \textbf{46} (1992) , 459--474. \newline
\newblock \href{http://www.ams.org/mathscinet-getitem?mr=1190349}{\texttt{\small MR1190349}},
          \href{http://dx.doi.org/10.1017/S0004972700012120}{\texttt{\small doi:10.1017/S0004972700012120}}

\bibitem{Ru94}
D.\ Ruelle,
\newblock \textit{Zeta Functions for Piecewise Monotone Maps of the Interval},
\newblock CRM Monograph Series, Vol.~4,
\newblock AMS, Providence, Rhode Island, 1994. \newline
\newblock \href{http://www.ams.org/mathscinet-getitem?mr=1274046}{\texttt{\small MR1274046}}

\bibitem{SS95}
T.A.\ Schmidt and M.\ Sheingorn,
\newblock \textit{Length spectra of the Hecke triangle groups},
\newblock Mathematische Zeitschrift \textbf{220} (1995), 369--397. \newline
\newblock \href{http://www.ams.org/mathscinet-getitem?mr=1362251}{\texttt{\small MR1362251}},
          \href{http://dx.doi.org/10.1007/BF02572621}{\texttt{\small doi:10.1007/BF02572621}}

\bibitem{Se89}
A.\ Selberg,
\newblock \textit{Remarks on the distribution of poles of Eisenstein series},
\newblock in \textit{Collected Papers}, Vol.\ 2, 15--46, Springer Verlag, 1989. \newline
\newblock \href{http://www.ams.org/mathscinet-getitem?mr=1295844}{\texttt{\small MR1295844}}

\bibitem{St08}
F.\ Str\"omberg,
\newblock \textit{Computation of Selberg's zeta functions on Hecke triangle groups}. \newline
\newblock \href{http://arxiv.org/abs/0804.4837}{\texttt{\small arXiv:0804.4837}}
\end{thebibliography}

\end{document}